\documentclass[11pt]{amsart} 
\usepackage[mathscr]{eucal}
\usepackage{amsmath,amsfonts}

\parskip=\smallskipamount

\newtheorem{theorem}{Theorem}[section]

\newtheorem{corollary}[theorem]{Corollary}

\newtheorem{remark}[theorem]{Remark}

\makeatletter
\@addtoreset{equation}{section}
\makeatother

\newcommand{\BB}{{\mathbb B}}
\newcommand{\CC}{{\mathbb C}}
\newcommand{\NN}{{\mathbb N}}

\newcommand{\ZZ}{{\mathbb Z}}

\newcommand{\RR}{{\mathbb R}}
\newcommand{\FF}{{\mathbb F}}
\newcommand{\TT}{{\mathbb T}}

\newcommand{\cA}{{\mathcal A}}

\newcommand{\cE}{{\mathcal E}}

\newcommand{\cG}{{\mathcal G}}
\newcommand{\cH}{{\mathcal H}}
\newcommand{\cK}{{\mathcal K}}

\newcommand{\cM}{{\mathcal M}}
\newcommand{\cN}{{\mathcal N}}
\newcommand{\cO}{{\mathcal O}}
\newcommand{\cP}{{\mathcal P}}
\newcommand{\cR}{{\mathcal R}}

\newdimen\expt
\def\boxit#1{\setbox0\hbox{$\displaystyle{#1}$}
      \hbox{\lower.4\expt
 \hbox{\lower3\expt\hbox{\lower\dp0
      \hbox{\vbox{\hrule height.4\expt
 \hbox{\vrule width.4\expt\hskip3\expt
      \vbox{\vskip3\expt\box0\vskip2\expt}%
 \hskip3\expt\vrule width.4\expt}\hrule height.4\expt}}}}}}
 
\begin{document}

 


\title [ Multivariable moment problems ] 
{ Multivariable moment problems   } 
 \author{Gelu Popescu}
\date{January 7, 2002 (revised version: September 11, 2002)}
\thanks{The author was
partially supported by an NSF grant}

\address{Department of Mathematics, The University of Texas 
at San Antonio \\ San Antonio, TX 78249, USA}
\email{\tt gpopescu@math.utsa.edu}

\begin{abstract}  In this paper we solve 
  moment problems 
for Poisson transforms and, more generally, for completely positive
  linear maps
   on unital ~$C^*$-algebras generated by ``universal'' row contractions
   associated with ~$\FF_n^+$, the free semigroup with $n$ generators.
   This class of ~$C^*$-algebras  includes
   the Cuntz-Toeplitz algebra  
~$C^*(S_1,\ldots, S_n)$  (resp.~$C^*(B_1,\ldots, B_n)$) 
generated by the creation operators on the full 
(resp.~symmetric, or anti-symmetric)) Fock
space with $n$ generators. 
As consequences, we  obtain   characterizations for  the orbits of contractive  
Hilbert modules over complex  free semigroup algebras 
 such as $\CC \FF_n^+$,\ $\CC[z_1,\ldots, z_n]$,
and, more generally, the quotient algebra $\CC \FF_n^+/J$,
 where $J$ is an arbitrary 
two-sided ideal of 
$\CC \FF_n^+$.

 All  these results are  extended to
 the generalized Cuntz algebra ~$\cO(*_{i=1}^n G_i^+)$, 
  where $G_i^+$ are the positive cones
 of discrete subgroups $G_i^+$ of the real line $\RR$.
 Moreover, we  characterize the orbits of Hilbert modules over  the quotient algebra
 $\CC *_{i=1}^n G_i^+/J$, where $J$ is an arbitrary two-sided ideal of the  
 free  semigroup algebra
$\CC *_{i=1}^n G_i^+$. 
\end{abstract}

\maketitle

\section{Introduction and Preliminaries}\label{INTRO}

Let $H_n$ be an $n$-dimensional complex  Hilbert space with orthonormal basis
$e_1,e_2,\dots,e_n$, where $n\in \{1,2,\dots\}$ or $n=\infty$.
  We consider the full Fock space  of $H_n$ defined by
$$F^2(H_n):=\bigoplus_{k\geq 0} H_n^{\otimes k},$$ 
where $H_n^{\otimes 0}:=\CC 1$ and $H_n^{\otimes k}$ is the (Hilbert)
tensor product of $k$ copies of $H_n$.
Define the left creation 
operators $S_i:F^2(H_n)\to F^2(H_n), \  i=1,\dots, n$,  by
$$
 S_i\psi:=e_i\otimes\psi, \  \psi\in F^2(H_n).
$$
The noncommutative analytic Toeplitz algebra $F_n^\infty$  (resp. 
noncommutative disc algebra $\cA_n$) is
 the
  weakly-closed (resp. the non-selfadjoint norm-closed) algebra generated by  
   $S_1,\dots, S_n$  
    and the identity.
 These algebras were introduced in 
\cite{Po-von} (see also \cite {Po-disc}) in connection with a multivariable
 noncommutative  von Neumann type inequality.
Let  $\FF_n^+$ be the free semigroup with $n$ generators $g_1,\dots, g_n$ and
  neutral element $g_0$. 
The length of $\alpha\in\FF_n^+$ is defined by
$|\alpha|:=k$, if $\alpha=g_{i_1}g_{i_2}\cdots g_{i_k}$, and
$|\alpha|:=0$, if $\alpha=g_0$.
We also define
$e_\alpha :=  e_{i_1}\otimes e_{i_2}\otimes \cdots \otimes e_{i_k}$  
 and $e_{g_0}= 1$.
It is  clear that 
$\{e_\alpha:\alpha\in\FF_n^+\}$ is an orthonormal basis of $F^2(H_n)$.
 If $T_1,\dots,T_n\in B(\cH)$  (the algebra of all bounded  
  linear operators on the Hilbert space $ \cH$), define 
$T_\alpha :=  T_{i_1}T_{i_2}\cdots T_{i_k}$,
if $\alpha=g_{i_1}g_{i_2}\cdots g_{i_k}$ and 
$T_{g_0}:=I$.

 Let $C^*(S_1,\dots,S_n)$ be the
$C^*$-algebra generated by $S_1,\dots,S_n$, the extension through 
compacts of the Cuntz algebra ${\cO}_n$ (see \cite{Cu}).
A map $\mu:C^*(S_1,\ldots, S_n)\to B(\cH)$ is called Poisson transform on the
Cuntz-Toeplitz algebra $C^*(S_1,\ldots, S_n)$   
if it is a unital completely positive
 linear  map such that 
  $$
  \mu(AX)=\mu(A)\mu(X),\quad \text{ for any } X\in \cA_n \cA_n^* 
 \text{ and } A\in \cA_n.
 $$
Notice that $\mu|\cA_n$ is a completely contractive representation of
the noncommutative disc algebra on the Hilbert space $\cH$.
It was shown in \cite{Po-poisson} that 
$\mu:C^*(S_1,\ldots, S_n)\to B(\cH)$ is a Poisson transform  if and only if 
there is a row contraction $T:=[T_1, \dots, T_n]$, i.e.,
$$
T_1T_1^*+\cdots + T_n T_n^*\leq I,
$$
  such that
$$
\mu(f)=\mu_T(f):=\lim_{r\to 1} K_r(T)^*[f\otimes I]K_r(T)
 $$
(in the uniform topology of $B(\cH)$), where
the Poisson Kernel  
$K_r(T):\cH\to F^2(H_n)\otimes\cH$, $0<r < 1$, is
 defined by
$$
K_r(T)h:= \sum\limits_{\gamma\in \FF^+_n} e_\gamma\otimes
\left(r^{|\gamma|}\Delta_r(T)T_\gamma^*h\right),
$$
and $\Delta_r(T):=(I-\sum_{i=1}^n r^2 T_iT_i^*)^{\frac{1} {2}}$.
In this case we have 
$$\mu_T(S_\alpha S_\beta^*)=T_\alpha T_\beta^*, \quad \alpha,\beta\in\FF_n^+. $$  
According to the noncommutative dilation theory
for row contractions (see \cite{Po-isometric}, \cite{Po-poisson}),
 $\mu_T|\cA_n^*$
is a subcoisometric representation on $\cA_n^*$. On the other hand, 
if $[T_1, \dots, T_n]$ is a  row isometry, then $\mu_T$ is a $*$-representation 
of the Cuntz-Toeplitz algebra $C^*(S_1,\ldots, S_n)$.

In the particular case when  $T_i T_j=T_j T_i$,  for any $ i,j=1,\ldots, n$,
 a similar  result  is true (see \cite{Po-poisson})
 if we replace the left creation operators $S_i$,  ~$i=1,\ldots, n$, 
 by their compressions $B_i:= P_{F^2_s} S_i|F^2_s$, \ $i=1,\ldots n$,
 to the symmetric Fock space $F^2_s\subset F^2(H_n)$. 
  In this commutative setting, 
  a  linear map $\rho:C^*(B_1,\ldots, B_n)\to B(\cH)$ 
  is a Poisson transform  if and only if 
there is a row contraction $T:=[T_1, \dots, T_n]$ with commuting entries,
  such that
  $$\rho(B_\alpha B_\beta^*)=T_\alpha T_\beta^*, \quad \alpha,\beta\in\FF_n^+.
   $$
  We refer to \cite{Po-poisson} for more information on noncommutative
   (resp.~commutative) Poisson transforms.
    The commutative case was extensively studied   by Arveson in \cite{Arv},
    where 
the creation operators $B_1,\ldots, B_n$ on the symmetric Fock space $F^2_s$
 are realized as multiplication operators by $z_1,\ldots, z_n$  (a system of
  coordinate functions for $\CC^n$) on the reproducing kernel Hilbert 
  space with reproducing kernel
  $K_n:\BB_n\times \BB_n\to \CC$ defined by
  $$
  K_n(z,w):= \frac{1}{1-\langle z,w\rangle}_{\CC^n},
  $$
  where $\BB_n$ is the open unit ball of $\CC^n$. 
  
 Recently,  the Poisson transforms have played a very important role
 in multivariable operator theory (see \cite{Po-poisson}, \cite{Po-tensor},
 \cite{Po-curvature}, 
 \cite{Po-similarity}, \cite{Arv},
  \cite{Arv2}, \cite{BB}, \cite{BBD}, \cite{Ar}).
  Moreover,   they  
  were essentially used in \cite{ArPo1} and  \cite{ArPo2} to prove
   interpolation results of Nevanlinna-Pick
     type for the unit ball of $\CC^n$.

Let $\CC \FF_n^+$ be the complex free semigroup algebra  generated
 by $\FF_n^+$ . 
Any $n$-tuple $T_1,\dots, T_n $ of bounded operators on a  Hilbert space
 $\cH$ gives rise to  a Hilbert  (left) module over  $\CC \FF_n^+$ in 
 the natural way
$$
f\cdot h:= f(T_1,\dots, T_n)h, \quad f\in  \CC \FF_n^+, h\in \cH,
$$
and $T:=[T_1,\dots, T_n] $ is a row contraction iff $\cH$ is a contractive 
$\CC \FF_n^+$-module, i.e.,
$$
\|g_1\cdot h_1+\cdots +g_n\cdot h_n\|^2\leq \|h_1\|^2+\cdots +\|h_n\|^2,
\quad
 h_1,\dots, h_n\in \cH.
$$
We   remark that each row contraction $T$  corresponds to a unique 
 contractive Hilbert modules over
 the free semigroup algebra $\CC \FF_n^+$, and 
to a unique Poisson transform $\mu_T$.
 A similar
 result holds true in the commutative case 
 when the contractive Hilbert
 modules are considered over $\CC [z_1,\ldots, z_n]$, the complex unital 
 algebra
  of all polynomials in $n$ commuting variables
 (see \cite{Po-poisson}, \cite{Po-curvature}, \cite{Arv}, \cite{Arv2}).
 Recently, 
  the curvature invariant for Hilbert modules over  $\CC [z_1,\ldots, z_n]$
  (resp. $\CC \FF_n^+$) was introduced and studied in \cite{Arv2} (resp. 
  \cite{Po-curvature} and \cite{Kr}). The Poisson  transforms played a 
  crucial
  role in connection with 
  numerical invariants associated with
  Hilbert modules  
 over  $\CC \FF_n^+$, 
  the complex free semigroup algebra 
   generated by the free semigroup on $n$ generators (see \cite{Po-curvature}
   and \cite{Po-similarity}).
We refer to \cite{Arv1}, \cite{P-book}, \cite{Pi-book} for basic facts concerning 
operator spaces and their completely positive (c.p.) (resp.~bounded) maps.

 In this paper we continue the study of Poisson transforms and, more generally,
 completely positive linear maps on  unital 
 ~$C^*$-algebras generated by ``universal'' row contractions, in connection with their
 moments.
 Our objective is to  obtain analogues of some classical moment problems 
 \cite{Ak}
 in our multivariable setting.
 
We  consider moment problems associated with any subset
 $\Sigma\subset \FF_n^+$ 
which is admissible in the sense that if 
$\alpha \beta\in \Sigma$, then $\beta\in \Sigma$.
A map $L:\Sigma \to B(\cH)$  (or its linear extension $L:\CC\Sigma \to B(\cH)$)
 is called
 {\it operator-valued moment map}  for Poisson transforms (resp.~c.p.~maps) on 
  $C^*(S_1,\ldots, S_n)$ if there is a Poisson transform (resp.~c.p.~map)
 $\mu:C^*(S_1,\ldots, S_n)\to B(\cH)$ such that
 $$
 L(p(g_1,\ldots, g_n))=\mu(p(S_1,\ldots, S_n))L(g_0)
 $$
 for any polynomial  $p(g_1,\ldots, g_n)$ in 
 $\CC\Sigma:=
 \{ \sum a_\sigma \sigma:\ a_\alpha \in \CC, \sigma\in \Sigma \}$.
 We call $\mu$ a representing Poisson transform (resp.~c.p.~map) for $L$.

In Section \ref{Moment},   we solve 
the operator-valued (resp.~vector-valued) moment problem 
for Poisson transforms  on  the Cuntz-Toeplitz algebra    
$C^*(S_1,\ldots, S_n)$.
As consequences, we obtain necessary and sufficient conditions
 for the existence
of a Poisson transform (resp.~$*$-representation) with
prescribed restrictions.
  We  also obtain a characterization for the orbits
  $
  \{f\cdot h: \ f\in  \FF_n^+\}$, \ $ h\in \cH$,
  of contractive Hilbert modules over $\CC \FF_n^+$, and show how to 
    construct contractive Hilbert
    modules with prescribed
   orbits.
   More precisely, we characterize the orbits
   $\{T_\sigma h:\ \sigma\in \FF_n^+\}$ generated by row contractions
   $[T_1,\ldots, T_n]$, ~$T_i\in B(\cH)$, and $h\in \cH$.
   On the other hand, we obtain a 
      characterization of 
   the orbits of   
 the $*$-representations of the $C^*$-algebra  $C^*(S_1,\ldots, S_n)$.

In Section \ref{commutative},  we solve 
the operator-valued (resp.~vector-valued) moment problem 
for Poisson transforms  on   $C^*$-algebras    
$C^*(B_1,\ldots, B_n)$ generated by universal row contractions
 $[B_1,\ldots, B_n]$ satisfying the commutation relations 
\begin{equation*}
  B_jB_i=\lambda_{ji} B_iB_j, \quad  1\leq i<j\leq n,
  \end{equation*}
where $\lambda_{ji}\in \CC$, \ $ 1\leq i<j\leq n$.
 We remark that if $\lambda_{ji}=1$  (resp.~$\lambda_{ji}=-1$),
 \ $1\leq i<j\leq n$, then $B_1,\ldots,
 B_n$ are the creation operators on the  symmetric (resp.~anti-symmetric)
  Fock space. 
 As a consequence (when $\lambda_{ji}=1$), 
 we  characterize the orbits
  $
  \{f\cdot h: \ f\in  \ZZ_n^+\}$, \ $ h\in \cH$,
  of contractive Hilbert modules over $\CC[z_1,\ldots, z_n]$.
On the other hand,
  the classical 
operatorial trigonometric moment problem  \cite{Ak} is extended 
to the commutative semigroup $\ZZ_n^+$.
More precisely, we solve the operator-valued trigonometric moment problem for 
completely positive linear maps on the $C^*$-algebras 
    $C^*(B_1,\ldots, B_n)$ described above (see Theorem \ref{trig-moment}).

In Section \ref{extension}, we present  moment problems for 
quotients of the free semigroup algebra
$\CC \FF_n^+$. Given a set $\cP_h$ of homogeneous polynomials in $\CC \FF_n^+$,
we show that there is a universal row contraction 
$[B_1,\ldots, B_n]$ such that $p(B_1,\ldots, B_n)=0$ for any
 $p\in \CC \FF_n^+$. Then  we solve the  moment problem for  Poisson 
 transforms  
   on the unital $C^*$-algebra  generated by 
 $B_1,\ldots, B_n$, and the identity.
 This provides a characterization 
   for  the orbits of contractive  
Hilbert modules over  the quotient algebra $\CC \FF_n^+/J_h$,
 where $J_h$ is  the 
two-sided ideal  generated by $\cP_h$.
On the other hand,  we solve 
 the operator-valued trigonometric moment problem for 
completely positive linear maps on $C^*(B_1,\ldots, B_n)$
 (see Theorem \ref{trig2}).
 
Finally, in Section \ref{moment-Cuntz}, we  show that all of our
 moment type results
can be extended to 
the generalized Cuntz algebra  $\cO(*_{i=1}^n G_i^+)$,
  where $G_i^+$ are the positive cones
 of discrete subgroups $G_i$ of the real line $\RR$. 
 Moreover, we  characterize the orbits of contractive Hilbert modules over
 the quotient algebra
 $\CC *_{i=1}^n G_i^+/J$, where $J$ is an arbitrary two-sided ideal of the  
 free  semigroup algebra
$\CC *_{i=1}^n G_i^+$.
We mentioned that
 the generalized Cuntz algebra $\cO(*_{i=1}^n G_i^+)$  coincides with 
the C*-algebra generated
by the left regular representation of the free proguct semigroup
 $*_{i=1}^n G_i^+$, when $n\geq 2$ and at least one subgroup $G_i$ 
 is dense in $\RR$ (see \cite{DPo}). 
   These algebras were studied by many authors (see for example \cite {Do},
    \cite{Cu},  
   \cite{D1}, 
   \cite{D2}, \cite{Ni},
    \cite{DPo}, \cite{LR}, and \cite{L}).

\bigskip

\section{Moment problems for the free semigroup $\FF_n^+$}
 \label{Moment}

 Let $\Sigma$ be an arbitrary set.
 An operator-valued function $K:\Sigma\times \Sigma\to B(\cH)$ is called 
 Hermitian kernel if $K(\tau,\sigma)=K(\sigma, \tau)^*$. 
  It is said to be positive
 semidefinite ($K\geq 0$) provided that 
 $$
 \sum_{\sigma, \tau\in \Sigma} \langle K(\sigma, \tau)h(\tau),
  h(\sigma)\rangle\geq 0
 $$
 for all finitely supported functions $h:\Sigma\to \cH$.
 If $K_1,K_2$ are operator-valued  kernels on $\Sigma$ we say that 
 $K_1\leq K_2$ if and only if $K_2-K_1$ is positive semidefinite.
 
In this section  
we  consider moment problems associated with any subset
 $\Sigma\subset \FF_n^+$ 
which is admissible in the sense that if 
$\alpha \beta\in \Sigma$, then $\beta\in \Sigma$.
A few examples that are of interest are the following:
\begin{enumerate}
\item[(i)]  $\Sigma= \FF_n^+$ is the full moment problem;
\item[(ii)]   $\Sigma=\{\sigma\in \FF_n^+:\ |\sigma|\leq m\}, ~m\in \ZZ^+$, 
is the truncation of order $m$;
\item[(iii)]   $\Sigma_\pi=\{\sigma\in \FF_n^+ : \ \varphi(\sigma)\in \Pi\}$,
 where
$\Pi$ is an admissible set for $\ZZ_n^+$ (see Section \ref {commutative}),
 is the abelian truncation;  
  \item[(iv)] 
   $\Sigma=\{ g_{i_1}g_{i_2}\cdots g_{i_k}, \ g_{i_2}\cdots g_{i_k},\ldots, 
   g_{i_k}, g_0\}$ is the truncation generated by
    $\omega:=g_{i_1}g_{i_2}\cdots g_{i_k}$,
    where $i_j\in \{1,\ldots, n\}$. 
\end{enumerate} 
A map $L:\Sigma \to B(\cH)$  (or its linear extension $L:\CC\Sigma \to B(\cH)$)
 is called
  operator-valued moment map  for
 Poisson transforms on 
  $C^*(S_1,\ldots, S_n)$ if there is a Poisson transform
 $\mu:C^*(S_1,\ldots, S_n)\to B(\cH)$ such that
 $$
 L(p(g_1,\ldots, g_n))=\mu(p(S_1,\ldots, S_n))L(g_0)
 $$
 for any polynomial  $p(g_1,\ldots, g_n)$ in 
 $\CC\Sigma:=
 \{ \sum a_\sigma \sigma:\ a_\alpha \in \CC,\  \sigma\in \Sigma \}$.
 We call $\mu$ a representing Poisson transform for $L$.
For each admissible set $\Sigma\subset \FF_n^+$,  we define 
\begin{equation}\label{lasi}
\Lambda_\Sigma:= \{(\alpha,\beta):\ \alpha, \beta \in \FF_n^+,
 \alpha\beta\in \Sigma\}.
\end{equation}
Notice that if $(\alpha,\beta)\in \Lambda_\Sigma$ and $\sigma\leq \alpha$, i.e.,
$\alpha=\sigma\tau$ for some $\tau\in \FF_n^+$, then 
$(\alpha\backslash \sigma,\beta)\in \Lambda_\Sigma$, where 
$\alpha\backslash \sigma:= \tau$. If $\tau\neq g_0$, denote $\sigma<\alpha$.
We associate with each map $L:\Sigma \to B(\cH)$ the Hermitian kernels
$K_j: \Lambda_\Sigma \times \Lambda_\Sigma \to B(\cH)$, ~$j=1,2$, defined by
\begin{equation}\label{k1}
K_1((\alpha,\beta), (\sigma,\gamma)):= 
L(\alpha\beta)^* L(\sigma\gamma) 
\end{equation}
and 
\begin{equation}\label{k2}
K_2((\alpha,\beta), (\sigma,\gamma)):= 
\begin{cases} L(\beta)^*L((\sigma \backslash \alpha)\gamma),
 & \text{ if } \alpha\leq \sigma\\
  L((\alpha\backslash \sigma)\beta)^* L(\gamma), 
    & \text{ if }  \sigma < \alpha \\
    0, & \text{ otherwise }.
\end{cases}
\end{equation}
Notice that  $K_1$ is a positive semidefinite kernel. The main result of this
section  is the following operator-valued moment problem for Poisson transforms
on the Cuntz-Toeplitz algebra  $C^*(S_1,\ldots, S_n)$.
  \begin{theorem}\label{main}
  Let $\Sigma\subseteq \FF_n^+$ be an admissible set and let
   $L:\Sigma\to B(\cH)$ be an operator-valued map. 
   Then $L$ is a  moment map for   Poisson transforms on the $C^*$-algebra 
    $C^*(S_1,\ldots, S_n)$ 
   if and only if $K_1\leq K_2$.
  \end{theorem} 
\begin{proof}
Assume that $\mu:C^*(S_1,\ldots, S_n)\to B(\cH)$ is a representing 
Poisson transform for $L$. Since $\mu|\cA_n$ is multiplicative, we have,
for any $(\alpha, \beta)\in \Lambda_\Sigma$, 
\begin{equation*}
L(\alpha \beta)=\mu(S_{\alpha \beta})L(g_0)=\mu(S_\alpha) L(\beta).
\end{equation*}
Let $[V_{1},\ldots,V_{n}]$ be the minimal isometric dilation of
 the row contraction $[\mu(S_1),\ldots, \mu(S_n)]$ on a Hilbert 
 space $\cK\supset \cH$ (see \cite{Po-isometric}).
If $f: \Lambda_\Sigma\to \cH$ is finitely supported, then we have
\begin{equation*}\begin{split}
\sum_{(\alpha,\beta), (\sigma,\gamma) \in \Lambda_\Sigma}
\langle
K_1((\alpha,\beta),& (\sigma,\gamma))  f(\sigma,\gamma),
 f(\alpha, \beta) \rangle\\
 &= 
  \|\sum_{(\alpha,\beta) \in \Lambda_\Sigma}  L(\alpha\beta)f(\alpha, \beta)\|^2
  \\
  &=
   \|\sum_{(\alpha,\beta) \in \Lambda_\Sigma} 
    \mu(S_\alpha)L(\beta)f(\alpha, \beta)\|^2\\
    &=
    \|P_\cH\sum_{(\alpha,\beta) \in \Lambda_\Sigma} 
    V_\alpha L(\beta)f(\alpha, \beta)\|^2\\
    &\leq 
    \|\sum_{(\alpha,\beta) \in \Lambda_\Sigma} 
    V_\alpha L(\beta)f(\alpha, \beta)\|^2.
\end{split}
\end{equation*}
Since $V_1,\ldots, V_n$ are isometries with orthogonal ranges and 
$P_\cH V_\omega|\cH= \mu(S_\omega)$, $\omega\in \FF_n^+$,  we infer that

\begin{equation*}\begin{split}
\sum_{(\alpha,\beta), (\sigma,\gamma) \in \Lambda_\Sigma}
  \langle V_\sigma L(\gamma)f(\sigma,\gamma), V_\alpha L(\beta)
  f(\alpha, \beta)\rangle 
  &=
  \sum_{\alpha \leq\sigma}
  \langle V_{\sigma\backslash \alpha} L(\gamma)f(\sigma,\gamma),
   L(\beta)f(\alpha, \beta)\rangle\\
  &+
  \sum_{\sigma<\alpha }
  \langle  L(\gamma)f(\sigma,\gamma), V_{\alpha\backslash \sigma}L(\beta)
  f(\alpha, \beta)\rangle\\
  &=
  \sum_{\alpha \leq\sigma}
  \langle \mu(S_{\sigma\backslash \alpha}) L(\gamma)f(\sigma,\gamma),
   L(\beta)f(\alpha, \beta)\rangle\\
  &+
  \sum_{\sigma<\alpha }
  \langle  L(\gamma)f(\sigma,\gamma), \mu(S_{\alpha\backslash \sigma})
  L(\beta)f(\alpha, \beta)\rangle\\
  &=
  \sum_{(\alpha,\beta), (\sigma,\gamma) \in \Lambda_\Sigma}
\langle
K_2((\alpha,\beta), (\sigma,\gamma))  f(\sigma,\gamma),
 f(\alpha, \beta) \rangle.
\end{split}
\end{equation*}
Therefore $K_1\leq K_2$ and the direct implication is proved. 

Conversely, assume $K_1\leq K_2$.
Let $\cH^{\Lambda_\Sigma}$ be the linear space of all finitely 
supported functions from $\Lambda_\Sigma$ to $\cH$.
Since $K_2\geq K_1\geq 0$, we can define a semidefinite form on 
$\cH^{\Lambda_\Sigma}$ by
$$
\langle f, g\rangle_{K_2}:=\sum_{(\alpha,\beta),
 (\sigma,\gamma)\in \Lambda_\Sigma} 
\langle
K_2((\alpha,\beta), (\sigma,\gamma))  f(\sigma,\gamma),
 f(\alpha, \beta) \rangle_\cH.
 $$
For each $i=1,\ldots, n$, define $V_i: \cH^{\Lambda_\Sigma}\to \cH^{\Lambda_\Sigma}$ by
$V_i f=\varphi$, where, for any $(\alpha, \beta)\in \Lambda_\Sigma$,
\begin{equation}\label{var}
\varphi(\alpha, \beta)=
\begin{cases}  f(\gamma, \beta),
 & \text{ if } \alpha=g_i\tau\\
    0, & \text{ otherwise }.
\end{cases}
\end{equation}
Define $\Lambda_0:= \{(\tau, \beta)\in \Lambda_\Sigma: g_i\tau\beta\in \Sigma\}$ and notice that
$$
K_2((g_i\tau,\beta), (g_i\tau', \gamma))=K_2((\tau,\beta), (\tau', \gamma)),
$$
if $(g_i\tau,\beta), (g_i\tau', \gamma)\in \Lambda_\Sigma$.
Notice that if $f$ is supported on $\Lambda_0$, then
\begin{equation*}\begin{split}
\langle V_i f, V_i f \rangle_{K_2} &=
\sum_{(\tau,\beta), (\tau',\gamma) \in \Lambda_0}
\langle K_2((g_i\tau,\beta), (g_i\tau', \gamma)) \varphi (g_i\tau', \gamma),
\varphi (g_i\tau,\beta) \rangle_\cH\\
&=
\sum_{(\tau,\beta), (\tau',\gamma) \in \Lambda_0}
\langle K_2((\tau,\beta), (\tau', \gamma)) f (\tau', \gamma),
f (\tau,\beta) \rangle_\cH\\
&= \langle  f,  f \rangle_{K_2}.
\end{split}
\end{equation*}
On the other hand, if $f$ is supported on $\Lambda_\Sigma\backslash \Lambda_0$,
then  $\langle V_i f, V_i f \rangle_{K_2}=0$.

Let $\cK$ denote the Hilbert space completion of $\cH^{\Lambda_\Sigma}/ \cN$, where
$\cN$ is the subspace of   null  vectors in the seminorm
$\langle \cdot,\cdot \rangle_{K_2}$.  Now, $V_i$ can be extended 
to a partial isometry on $\cK$, which we denote also by $V_i$.
Notice that the range of $V_i$ is contained in the closed span of functions supported on
$$
\{(\alpha, \beta)\in \Lambda_\Sigma:\ \alpha=g_i\tau \text{ for some } 
\tau\in \FF_n^+\}.
$$
According to the definition of $K_2$, we have
$$
K_2((g_i\tau,\beta), (g_j\tau', \gamma))=0 \quad \text{ if } i\neq j.
$$
Hence, the partial isometries $V_i$, \ $i=1,\ldots, n$,
 have orthogonal ranges and therefore,
 \begin{equation}\label{rowiso}
 V_1V_1^*+\cdots + V_n V_n^*\leq I.
 \end{equation}
Define $X:\cH^{\Lambda_\Sigma}\to \cH$ by
\begin{equation}\label{xf}
Xf:= \sum_{(\alpha, \beta)\in \Lambda_\Sigma} L(\alpha \beta)f(\alpha, \beta),
 \quad
 f\in \cH^{\Lambda_\Sigma}.
\end{equation}
Notice that, since $K_1\leq K_2$, we have
\begin{equation*}
\begin{split}
\|Xf\|_\cH^2&= \|\sum_{(\alpha, \beta)\in \Lambda_\Sigma}
 L(\alpha \beta)f(\alpha, \beta)\|_\cH^2\\
&=\sum_{(\alpha, \beta), (\sigma,\gamma)\in \Lambda_\Sigma}
\langle
K_1((\alpha,\beta), (\sigma,\gamma))  f(\sigma,\gamma),
 f(\alpha, \beta) \rangle_\cH \\
 &\leq\sum_{(\alpha, \beta), (\sigma,\gamma)\in \Lambda_\Sigma}
\langle
K_2((\alpha,\beta), (\sigma,\gamma))  f(\sigma,\gamma),
 f(\alpha, \beta) \rangle_\cH \\
 &=\langle f, f \rangle_{K_2}.
 \end{split}
\end{equation*}
Therefore, $X$ extends to a contraction from $\cK$ to $\cH$, 
which we denote also by $X$.
We need to show that, if $\sigma\in \Sigma, h\in \cH$, then 
\begin{equation}\label{x*l}
X^*L(\sigma)h= \psi_{(g_0, \sigma), h},
\end{equation}
where 
\begin{equation*}
\psi_{(g_0, \sigma), h}(\alpha, \beta):=
\begin{cases}  h,
 & \text{ if } \alpha=g_0, ~\beta=\sigma\\
    0, & \text{ otherwise }.
\end{cases}
\end{equation*}
Indeed, for any $f\in \cH^{\Lambda_\Sigma}$, we have
\begin{equation*}
\begin{split}
\langle f, X^*L(\sigma)h \rangle_{K_2}&=\langle Xf, L(\sigma)h \rangle_\cH\\
&=
\sum_{(\omega, \eta) \in \Lambda_\Sigma}
\langle L(\omega \eta)  f(\omega, \eta), L(\sigma)h\rangle_\cH\\
&=\sum_{(\omega, \eta) \in \Lambda_\Sigma}
\langle K_2((g_0,\sigma), (\omega,\eta))  f(\omega,\eta), h \rangle_\cH\\
&=
\sum_{(\alpha, \beta), (\omega, \eta) \in \Lambda_\Sigma}
\langle K_2((\alpha,\beta), (\omega,\eta))  f(\omega,\eta),
 \psi_{(g_0, \sigma), h}(\alpha, \beta) \rangle_\cH\\
 &= \langle f, \psi_{(g_0, \sigma), h}\rangle_{K_2}.
\end{split}
\end{equation*}
Now, define $T_i:= XV_i X^*\in B(\cH)$, \ $i=1,\ldots, n$. 
Since $X$ is a contraction
 and relation
\eqref{rowiso} holds we infer that
\begin{equation}\label{row}
 T_1T_1^*+\cdots + T_n T_n^*\leq I.
 \end{equation}
Using relation \eqref{x*l}, we deduce that 
for any $\sigma\in \Sigma, h\in \cH$,  
\begin{equation}\label{vx*l}
V_iX^*L(\sigma)h= \psi_{(g_i, \sigma), h},
\end{equation}
where 
\begin{equation*}
\psi_{(g_i, \sigma), h}(\alpha, \beta):=
\begin{cases}  h,
 & \text{ if } \alpha=g_i, ~\beta=\sigma\\
    0, & \text{ otherwise }.
\end{cases}
\end{equation*}
By relations \eqref{x*l} and \eqref{xf}, we obtain
 \begin{equation*}
 T_i L(\sigma) h= 
  \begin{cases}
 L(g_i\sigma)h, &  \text{ if } g_i\sigma\in \Sigma\\
  0, & \text{ otherwise }.
\end{cases}
\end{equation*}
 Now, it is clear that this relation implies
\begin{equation}\label{tomega}
T_\omega L(g_0) h=T_{i_1}\cdots T_{i_k}L(g_0) h=L( g_{i_1}\cdots g_{i_k})h,
\end{equation}
for any $\omega:= g_{i_1}\cdots g_{i_k}\in \Sigma$ and $h\in \cH$.
Hence, $T_\omega L(g_0)=L(\omega)$ for any $\omega\in \Sigma$.

Let $\mu:C^*(S_1,\ldots, S_n)\to B(\cH)$ be the Poisson transform associated with 
the row contraction 
$[T_1,\ldots, T_n]$. Since $\mu(S_\alpha)= T_\alpha$,\ $\alpha\in \FF_n^+$,
 relation \eqref{tomega}
implies
$$
 L(p(g_1,\ldots, g_n))=\mu(p(S_1,\ldots, S_n))L(g_0)
 $$
 for any polynomial  $p(g_1,\ldots, g_n)$ in 
 $\CC\Sigma$. The proof is complete.
\end{proof}

  \begin{corollary}\label{vect}
  Let $\cM$ be a subspace of  a Hilbert space  $\cH$ and let
   $\Sigma\subseteq \FF_n^+$ be an admissible set.  
    Given a map  $\Gamma:\Sigma\to B(\cM,\cH)$, there is a Poisson transform 
    $\mu: C^*(S_1,\ldots, S_n)\to B(\cH)$ such that
    $$
    \Gamma(p(g_1,\ldots, g_n))=\mu(p(S_1,\ldots, S_n)\Gamma(g_0)
    $$
    for any polynomial  $p(g_1,\ldots, g_n)$ in 
 $\CC\Sigma$, 
   if and only if $K_1\leq K_2$, where the kernels $K_1, K_2$ are 
   associated with
   the operator-valued map $L:\Sigma\to B(\cH)$ defined
   by 
   $$L(\sigma):= \Gamma(\sigma) P_\cM, \quad \sigma\in \FF_n^+,
   $$
   where $P_\cM$ is the orthogonal projection from  $\cH$ onto $\cM$.
  \end{corollary} 

In particular, if $\Gamma(g_0)$ is the inclusion of $\cM$ into $\cH$,
 then we obtain necessary
and sufficient conditions for the existence of a Poisson transform 
 $\mu: C^*(S_1,\ldots, S_n)\to B(\cH)$ with prescribed restriction 
 to  $\cM$, i.e., 
 $$
    \Gamma(p(g_1,\ldots, g_n))=\mu(p(S_1,\ldots, S_n)|\cM
    $$
    for any polynomial  $p(g_1,\ldots, g_n)$ in 
 $\CC\Sigma$.
 
Now, we can obtain a vector-valued version of  
Theorem \ref{main}.
A map $M:\Sigma \to \cH$  (or its linear extension $M:\CC\Sigma \to \cH$)
 is called
 {\it vector-valued moment map}  for Poisson transforms on $C^*(S_1,\ldots, S_n)$ if there is a Poisson transform
 $\mu:C^*(S_1,\ldots, S_n)\to B(\cH)$ such that
 $$
 M(p(g_1,\ldots, g_n))=\mu(p(S_1,\ldots, S_n))M(g_0)
 $$
 for any polynomial  $p(g_1,\ldots, g_n)$ in 
 $\CC\Sigma$.
 We call $\mu$ a representing Poisson transform for $M$.
We associate with each map $M:\Sigma \to \cH$ the hermitian kernels
$K_j: \Lambda_\Sigma \times \Lambda_\Sigma \to \CC$, ~$j=1,2$, defined by
\begin{equation}\label{k3}
K_3((\alpha,\beta), (\sigma,\gamma)):= 
\left< M(\sigma\gamma), M(\alpha\beta)\right>
\end{equation}
and 
\begin{equation}\label{k4}
K_4((\alpha,\beta), (\sigma,\gamma)):= 
\begin{cases} \left< M((\sigma \backslash \alpha)\gamma), M(\beta)\right>,
 & \text{ if } \alpha\leq \sigma\\
  \left< M(\gamma), M((\alpha\backslash \sigma)\beta)\right>, 
    & \text{ if }  \sigma < \alpha \\
    0, & \text{ otherwise}.
\end{cases}
\end{equation}

  \begin{theorem}\label{main-vect}
  Let $\Sigma\subseteq \FF_n^+$ be an admissible set and let
   $M:\Sigma\to \cH$ be
   a map with values in a Hilbert space. 
   Then $M$ is a vector-valued moment map for Poisson transforms on
    $C^*(S_1,\ldots, S_n)$ 
   if and only if $K_3\leq K_4$.
  \end{theorem} 
\begin{proof} In Corollary \ref{vect}, take $\cM:= \text{\rm span} \{ M(g_0)\}$ and 
let $\Gamma(g_0)$ be the inclusion of $\cM$ into $\cH$.   
The result follows 
noticing that, in this particular case,  $K_1\leq K_2$ if and only if $K_3\leq K_4$.
\end{proof}
As mention in the introduction,  in the particular case when 
$\Sigma=\FF_n^+$, Theorem \ref{main-vect} provides
 a characterization for the orbits
$
\{f\cdot h:\ f\in \FF_n^+\}, \quad h\in \cH,
$
of contractive Hilbert modules over $\CC \FF_n^+$.
In what follows, we will characterize the orbits 
 of $*$-representations of the Cuntz-Toeplitz algebra $C^*(S_1,\ldots, S_n)$.
  \begin{theorem}\label{rep}
  Let $\Sigma\subseteq \FF_n^+$ be an admissible set and let
   $L:\Sigma\to B(\cH)$ be an operator-valued map. 
   Then  there is a $*$-representation $\pi: C^*(S_1,\ldots, S_n) \to B(\cE)$, 
   $\cE\supset \cH$, such that
   \begin{equation}\label{pi}
    L(p(g_1,\ldots, g_n))=\pi(p(S_1,\ldots, S_n)L(g_0)
    \end{equation}
    for any polynomial  $p(g_1,\ldots, g_n)$ in 
 $\CC\Sigma$, 
   if and only if $K_1= K_2$.
  \end{theorem} 
\begin{proof}
If $\pi:C^*(S_1,\ldots, S_n)\to B(\cE)$ is a a $*$-representation
 on a Hilbert space $\cE\supset \cH$,
 then 
$\pi(g_i)$, \ $i=1,\ldots, n$, are isometries with orthogonal ranges.
Then, for any $(\alpha,\beta), (\sigma,\gamma) \in \Lambda_\Sigma$ and $h, h'\in \cH$,
 relation 
\eqref{pi}  implies
\begin{equation*}\begin{split}
\langle
K_1((\alpha,\beta), (\sigma,\gamma))  h,
 h'\rangle 
 &= 
    \langle L(\sigma \gamma)h,  L(\alpha\beta)h'\rangle 
  \\
  &=
     \langle \pi(S_\sigma) L(\gamma)h,
    \pi(S_\alpha) L(\beta) h'\rangle\\
  &=\begin{cases}
  \langle \pi(S_{\sigma\backslash \alpha}) L(\gamma) h,
   L(\beta) h'\rangle, & \text{ if } \alpha\leq \sigma\\
  \langle  L(\gamma) h, \pi(S_{\alpha\backslash \sigma})
  L(\beta) h'\rangle, & \text{ if } \sigma <\alpha\\
  0, & \text{ otherwise }
  \end{cases}
  \\
 &=
\langle
K_2((\alpha,\beta), (\sigma,\gamma))   h,
 h' \rangle.
\end{split}
\end{equation*}
Therefore $K_1= K_2$ and the direct implication is proved. 

Conversely, assume $K_1= K_2$. Following the proof of Theorem \ref{main}, 
it is easy to see that, in this case, the operator $X:\cK\to \cH$ is an isometry
 and $T_i:=XV_iX^*$, \ $i=1,\ldots, n$,   are now partial isometries
  with orthogonal ranges.
 Moreover, the initial space of $V_i$ is the closed span of functions supported on 
 $$
 \Lambda_0:=\{ (\tau, \beta)\in \Lambda_\Sigma: \ g_i\tau\beta\in \Sigma\}.
 $$
 Taking into account the definition of the isometry $X$ (see \eqref{xf}), 
 we deduce that
  the initial space of $T_i$ is
 $$
 \text{\rm span} \{L(\sigma)h:\ g_i\sigma \in \Sigma, \ h\in \cH\}.
 $$   
Let $[W_1,\ldots, W_n]$ be the minimal isometric dilation of 
$[T_1,\ldots, T_n]$ on a Hilbert space $\cE\supset \cH$. 
Since $T_i$ is a partial isometry,
it is clear that
$W_ix=T_ix$  when $x$ is in the initial space of $T_i$. Hence
\begin{equation}\label{witi}
W_i L(\sigma) h=T_i L(\sigma) h
\end{equation}
 if $g_i\sigma\in \Sigma$ and $h\in \cH$.
Let $\pi:C^*(S_1,\ldots, S_n)\to B(\cE)$ be  the  $*$-representation
 generated by 
$\pi(S_i)=W_i$, \ $i=1,\ldots, n$.
Fix an arbitrary element $\omega:=g_{i_1}\ldots g_{i_k}$  in $\Sigma$.
 Since $\Sigma$ 
is admissible, the elements 
$g_{i_2}\ldots g_{i_k}, \ g_{i_3}\ldots g_{i_k}, \ldots, g_{i_k}$, and $ g_0$ 
are in $\Sigma$.
A repeated application of  the relations  \eqref{witi}  and \eqref{tomega} shows that
\begin{equation*}\begin{split}
W_\omega L(g_0)h&= W_{i_1}\cdots W_{i_k} L(g_0)h\\
&=W_{i_1}\cdots W_{i_{k-1}} T_{i_k} L(g_0)h
=W_{i_1}\cdots W_{i_{k-1}}  L(g_{i_k})h\\
&\ \cdots\cdots\\
&=L(g_{i_1}\cdots   g_{i_k})h.
\end{split}
\end{equation*}
Therefore, $\pi(S_\omega)L(g_0)= L(\omega)$ for any $\omega\in \Sigma$, which
completes the proof.  
\end{proof}

Given $L:\Sigma\to B(\cH)$, where $\Sigma$ is an admissible set in $\FF_n^+$,
 define the Hermitian kernels $K_1', K_2':\Sigma\times \Sigma\to B(\cH)$ by 
 setting
 $K_1'(\sigma, \tau):=L(\sigma)^* L(\tau)$ and 
  \begin{equation*}
  K_2'(\sigma, \tau):=
  \begin{cases}
  L(\tau\backslash \sigma), & \text{ if } \sigma\leq \tau\\
  L(\sigma\backslash \tau)^*, & \text{ if }\tau < \sigma\\
  0, & \text{ otherwise}.
  \end{cases}
  \end{equation*}

\begin{remark}\label{K'}  $K_1=K_2$ if and only if $K_1'=K_2'$.
\end{remark}
\begin{proof}
If $K_1=K_2$, then by setting $\beta=\gamma=g_0$ in relations \eqref{k1}
and \eqref{k2}, we infer  $K_1'=K_2'$.
Conversely, assume  $K_1'=K_2'$. Using the observation
 that $\omega\sigma\leq \omega \alpha$ if and only if 
 $\sigma\leq \alpha$, we infer that  
$
K_1'(\omega\sigma, \omega \alpha)= K_1'(\sigma, \alpha)
$
for any $\omega\sigma, \omega \alpha\in \Sigma$, and 
$
K_1'(g_i\sigma, g_j \alpha)=0
$
if $i\neq j$ and $g_i\sigma, g_j \alpha\in \Sigma$.
Now, using these relations, it is easy to see that, for any 
$(\alpha, \beta), (\sigma, \gamma)\in \Lambda_\Sigma$, 
\begin{equation*}
\begin{split}
K_1((\alpha, \beta), (\sigma, \gamma))&= 
K_1'(\alpha\beta, \sigma \gamma)\\
&=
\begin{cases}
  K_1'(\beta,(\sigma\backslash \alpha)\gamma), & \text{ if } \alpha\leq \sigma\\
   K_1'( (\alpha\backslash \sigma)\beta, \gamma), & \text{ if }\sigma < \alpha\\
  0, & \text{ otherwise}
  \end{cases}\\
  &=K_2((\alpha, \beta), (\sigma, \gamma)),
\end{split}
\end{equation*}
which completes the proof.
\end{proof}
Consequently, one can reformulate Theorem \ref{rep} in terms of the kernels
$K_1'$ and $K_2'$.
 We should mention that  one can  also get 
  an analogue of Corollary \ref{vect}   for
 $*$-representations of $C^*(S_1,\ldots, S_n)$.
Moreover, as a particular case, we obtain the following vector-valued version 
of Theorem \ref{rep}.

  \begin{theorem}\label{rep1}
  Let $\Sigma\subseteq \FF_n^+$ be an admissible set and let
   $M:\Sigma\to \cH$ be a map. 
   Then  there is a $*$-representation $\pi: C^*(S_1,\ldots, S_n) \to B(\cE)$, 
   $\cE\supset \cH$, such that
   \begin{equation*} 
    M(p(g_1,\ldots, g_n))=\pi(p(S_1,\ldots, S_n)M(g_0)
    \end{equation*}
    for any polynomial  $p(g_1,\ldots, g_n)$ in 
 $\CC\Sigma$, 
   if and only if $K_3= K_4$.
  \end{theorem} 
We associate with each map $M:\Sigma \to \cH$ the hermitian kernels
$K_j': \Sigma \times \Sigma \to \CC$, ~$j=3,4$, defined by
\begin{equation}\label{k3'}
K_3'((\alpha,\sigma):= 
\left< M(\sigma), M(\alpha)\right>
\end{equation}
and 
\begin{equation}\label{k4'}
K_4'((\alpha,\sigma):= 
\begin{cases} \left< M(\sigma \backslash \alpha), M(g_0)\right>,
 & \text{ if } \alpha\leq \sigma\\
  \left< M(g_0), M(\alpha\backslash \sigma)\right>, 
    & \text{ if }  \sigma < \alpha \\
    0, & \text{ otherwise}.
\end{cases}
\end{equation}
  One can easily       prove a result similar to
  Remark \ref{K'} for the kernels $K_3', K_4'$, and reformulate 
  Theorem \ref{rep1} in terms of these
   new kernels.
  We leave this task to the reader.

\bigskip

\section{Moment problems for the commutative semigroup $\ZZ_n^+$}
 \label{commutative}

 Let $\ZZ^+$ be the set of nonnegative integers and  
 $$
 \ZZ_n^+:=\{{\bf k}=(k_1,\ldots, k_n):\    k_j\in \ZZ^+, j=1,\ldots, n\}.
 $$
 If ${\bf k}=(k_1,\ldots, k_n)\in \ZZ^+_n$, we denote
  $|{\bf k}|:= k_1+\cdots + k_n$, and 
 if 
 $T:=(T_1,\ldots, T_n)$, \ $T_i\in B(\cH)$,  
  is an $n$-tuple of 
   operators, then we  
 set ~$T^{\bf k}:=T_{1}^{k_1}\cdots T_{n}^{k_n}$.
 Let $\CC[z_1,\ldots, z_n]$ be the complex unital algebra 
 of all polynomials in $n$ commuting variables. For any $\Pi\subset \ZZ_n^+$ ,
 denote 
 $$
 \CC[\Pi]:=\{p\in \CC[z_1,\ldots, z_n]:\ p=\sum_{{\bf k}\in \Pi} a_{\bf k}
  z^{\bf k}\},
 $$
 where $z^{\bf k}:= z_1^{k_1}\cdots z_n^{k_n}$. Identifying each ${\bf k}\in \Pi$
  with $z^{\bf k}$ we can consider $\Pi$ as a subset of $\CC [\Pi]$.
  The set $\ZZ_n^+$ is endowed with the product order $\ll$.
 We call a set  $\Pi\subset \ZZ_n^+$   admissible if ${\bf k}\in \Pi$,
 ${\bf m}\in \ZZ_n^+$, and $ ~{\bf m}\ll {\bf k}$, imply ${\bf m}\in \Pi$.
A few examples that are of interest are the following:

\begin{enumerate}
\item[(i)]  $\Pi= \ZZ_n^+$ is the full moment problem;
\item[(ii)]   $\Pi=\{{\bf k}\in \ZZ_n^+:\ |{\bf k}|\leq m\}, ~m\in \ZZ^+$, 
is the truncation of order $m$;
\item[(iii)]   $\Pi=\{{\bf k} \in \ZZ_n^+: \ {\bf k}\ll {\bf m}\}$
 is the truncation of order
  ${\bf m} \in \ZZ_n^+$.
  \item[(iv)] 
   $\Pi=\{{ \bf k}=(k_1,\ldots, k_n)\in \ZZ_n^+:\ 
    \text{\rm card}\{j:\ k_j\neq 0\}\leq p\}$,\  $p\in \ZZ^+$.
\end{enumerate} 
A map $\Gamma:\Pi \to B(\cH)$  (or its linear extension
 $\Gamma:\CC[\Pi] \to B(\cH)$)
 is called
 {\it operator-valued moment map}  for Poisson transforms on
  the Toeplitz 
 algebra $C^*(B_1,\ldots, B_n)$ if there is a Poisson transform
 $\rho:C^*(B_1,\ldots, B_n)\to B(\cH)$ such that
 \begin{equation}\label{rho}
 \Gamma(\sum a_{\bf k} z^{\bf  k})=\rho(\sum a_{\bf k} B^{\bf k})\Gamma({\bf 0})
 \end{equation}
 for any polynomial  $\sum a_{\bf k}  z^{\bf k}$ in 
 $\CC[\Pi]$, where $B:=(B_1,\ldots, B_n)$.
 We call $\rho$ a representing Poisson transform for $\Gamma$.

  There is a canonical homomorphism
$\varphi$ of $\FF_n^+=\ZZ^+\ast\cdots \ast\ZZ^+$ onto $\ZZ_n^+$ such that 
it is the identity on  each $\ZZ^+$.
Notice that if $\Pi$ is an admissible set for $\ZZ_n^+$, then
\begin{equation*}
\Sigma_\pi:= \{\sigma\in \FF_n^+:\ \varphi(\sigma)\in \Pi\}
\end{equation*}
is an admissible set for $\FF_n^+$.  
We associate with each map $\Gamma:\Pi \to B(\cH)$ 
a map $L_\Gamma:\Sigma_\pi\to B(\cH)$ by setting $L_\Gamma(\sigma):= 
\Gamma (\varphi(\sigma))$, 
$\sigma\in \Sigma_\pi$.

 \begin{theorem}\label{main-comm} Let $\Pi$ be an admissible set 
 in $\ZZ^+_n$ and let
 $\Gamma:\Pi \to B(\cH)$ be an operator-valued  map. Then $\Gamma$ is a 
  moment map  for  Poisson transforms on 
 the Toeplitz 
 algebra $C^*(B_1,\ldots, B_n)$  if and only if $K_1\leq K_2$, where the kernels
 $K_1$ and $K_2$ are associated with the admissible  set
  $\Sigma_\pi\subset \FF_n^+$ and the map $L_\Gamma:\Sigma_\pi\to B(\cH)$, 
  as defined by \eqref{k1} and \eqref{k2}.
 \end{theorem}
 \begin{proof}
 Assume $\rho: C^*(B_1,\ldots, B_n)\to B(\cH)$  is a representing Poisson 
 transform for $\Gamma$. According to \eqref{rho}, for any 
 $(\alpha, \beta)\in \Lambda_{\Sigma_\pi}$, we have
 \begin{equation}\label{lla}\begin{split}
 L_\Gamma(\alpha \beta)&= \Gamma(\varphi(\alpha)+\varphi(\beta))\\
 &= \rho(B^{\varphi(\alpha)+\varphi(\beta)})\Gamma({\bf 0})\\
 &= \rho(B^{\varphi(\alpha)}) \rho(B^{\varphi(\beta)})\Gamma({\bf 0})\\
 &= \rho(B^{\varphi(\alpha)}) \Gamma(\varphi(\beta))\\
 &= \mu(S_\alpha) L_\Gamma(\beta),
 \end{split}
 \end{equation}
 where $\mu: C^*(S_1,\ldots, S_n)\to B(\cH)$  is  the noncommutative  Poisson 
 transform associated with the row contraction $[\rho(B_1), \ldots, \rho(B_n)]$,
 i.e.,
 \begin{equation}\label{miu}
 \mu(S_\alpha S_\beta^*)= \rho(B^{\varphi(\alpha)}{B^{\varphi(\beta)}}^*),
  \quad \alpha, \beta \in \FF_n^+.
 \end{equation}
 Therefore, $\mu$ is a representing Poisson transform for $L_\Gamma$.
  According to Theorem \ref{main},  we infer that $K_1\leq K_2$.
 
 Conversely, assume that $K_1\leq K_2$. Following the proof of 
 Theorem \ref{main}, we find  the operators $T_i\in B(\cH)$, \ $i=1,\ldots, n$,
  defined by $T_i:= XV_iX^*$ with property that
  $$T_1T_1^*+\cdots +T_nT_n^*\leq I.
  $$
 Due to the definition of $X$ (see \eqref{xf}), its range is 
 included in the subspace
 $$
 \cM:=\text{\rm span} \{L_\Gamma(\sigma)h:\ \sigma\in \Sigma_\pi,
  ~h\in \cH\}.
 $$
 Hence, we have
 \begin{equation}\label{tim}
 T_i|\cM^\perp= 0, \quad i=1,\ldots, n.
 \end{equation}
 On the other hand, according to the definition of $T_i$, if 
 $\sigma\in \Sigma_\pi$
 and  ~$h\in \cH$, then
 \begin{equation*}
 T_iT_jL_\Gamma(\sigma)h=
 \begin{cases}   L_\Gamma(g_i g_j\sigma)h,
 & \text{ if }  g_ig_j\sigma\in \Sigma_\pi\\
    0, & \text{ otherwise}.
\end{cases}
 \end{equation*}
 Similarly, we have
 \begin{equation*}
 T_jT_iL_\Gamma(\sigma)h=
 \begin{cases}   L_\Gamma(g_j g_i\sigma)h,
 & \text{ if }  g_jg_i\sigma\in \Sigma_\pi\\
    0, & \text{ otherwise}.
\end{cases}
 \end{equation*}
 Due to the definition of $\Sigma_\pi$,  $ g_ig_j\sigma\in \Sigma_\pi$ if
  and only if 
 $ g_jg_i\sigma\in \Sigma_\pi$. On the other hand, according to the
  definition of $L_\Gamma$, we have 
$ L_\Gamma(g_i g_j\sigma)h= L_\Gamma(g_j g_i\sigma)h$, and therefore
 $T_iT_j|\cM=T_jT_i|\cM$. Now,  using  relation \eqref{tim}, we get 
 $$
 T_iT_j=T_jT_i, \quad i,j=1,\ldots, n.
 $$
 Since $[T_1,\ldots, T_n]$ is a row contraction with commuting entries, 
 according to Theorem 9.2 from
 \cite{Po-poisson}, there is a  unique Poisson transform
  $\rho: C^*(B_1,\ldots, B_n)\to B(\cH)$ such that
  \begin{equation}\label{robt}
  \rho(B^{\bf k} {B^{\bf m}}^*)=T^{\bf k} {T^{\bf m}}^*, \quad  
  {\bf k},  {\bf m}\in \ZZ_n^+,
  \end{equation}
  where $T:=(T_1,\ldots, T_n)$.
  For each ${\bf k}=(k_1,\ldots, k_n)\in \Pi$, 
  let $\omega_{\bf k}:= g_1^{k_1}\cdots g_n^{k_n}\in \FF_n^+$
  and notice that $\varphi(\omega_{\bf k})={\bf k}$.
  Moreover,  if ${\bf k}\in \Pi$, then 
  $\omega_{\bf k}\in \Sigma_\pi$ and, according to \eqref{tomega}, we have
  \begin{equation}\label{lla1}
  L_\Gamma(\omega_{\bf k})=\mu(S_{\omega_{\bf k}})L_\Gamma(g_0),
  \end{equation}
  where  $\mu: C^*(S_1,\ldots, S_n)\to B(\cH)$  is 
   the noncommutative  Poisson 
 transform associated with the row contraction $[T_1,\ldots, T_n]$. 
   Using  the definition of $L_\Gamma$ and relations \eqref{robt}, 
     \eqref{lla1},
    we obtain
    \begin{equation*}\begin{split}
    \rho(B^{\bf k}) \Gamma({\bf 0})&= T_1^{k_1}\cdots T_n^{k_n} L_\Gamma(g_0)\\
    &=\mu(S_{\omega_{\bf k}})L_\Gamma(g_0)=L_\Gamma(\omega_{\bf k})\\
    &=\Gamma({\bf k}),
    \end{split}
    \end{equation*}
for any ${\bf k}\in \Pi$. The proof is complete.
 \end{proof}
 Let us remark that
  an analogue of Corollary
 \ref{vect} holds also   in this  commutative case.
 We associate with any map $N:\Pi \to \cH$ 
a map $M_N:\Sigma_\pi\to \cH$ by setting $M_N(\sigma):= 
N (\varphi(\sigma))$, 
$\sigma\in \Sigma_\pi$.
 One can obtain a vector-valued version of Theorem \ref{main-comm}.
 The proof is similar to that of Theorem \ref{main-vect}  and
 uses  Theorem \ref{main-comm}. We should omit it.
 \begin{theorem}\label{vecom} Let $\Pi$ be an admissible set in $\ZZ^+_n$ and
  let
 $N:\Pi \to \cH$ be a map with values in a Hilbert space. Then $N$ is a 
 vector-valued moment map  for Poisson transforms on  the Toeplitz 
 algebra $C^*(B_1,\ldots, B_n)$  if and only if $K_3\leq K_4$, where the kernels
 $K_3$ and $K_4$ are associated with the admissible  set
  $\Sigma_\pi\subset \FF_n^+$ and the map $M_N:\Sigma_\pi\to \cH$, 
  as defined by \eqref{k3} and \eqref{k4}.
 \end{theorem}
 Notice that, in the particular case when $\Pi=\ZZ_n^+$,  Theorem \ref{vecom} provides a
   characterization for  the orbits
  $
  \{f\cdot h: \ f\in  \ZZ_n^+\}$, \  $h\in \cH$,
  of contractive Hilbert modules over $\CC[z_1,\ldots, z_n]$.

 A $B(\cH)$-valued semispectral measure on 
 $\TT:=\{z\in \CC:\ |z|=1\}$ is a
  linear positive map from $C(\TT)$, the set of
   continous functions on the unit circle, into $B(\cH)$.
    Since  $C(\TT)$ is commutative $\mu$ is completely positive.
The truncated  operatorial trigonometric moment problem asks if, 
given the operators $A_k\in B(\cH)$, $ k=0,1,\dots,m$,  $(A_0=I)$, 
there exists a semispectral measure on
 $\TT$ such that $A_k=\mu(e^{ikt})$,\ $k=0,1,\dots,m$. 
 The answer is provided in
terms of the positivity of an associated  Toeplitz operator.
 See \cite{Ak} for more information on the classical moment problem. 
 A non-commutative analogue of this problem was obtained in 
 \cite{Po-positive} for the Cuntz-Toeplitz algebra $C^*(S_1,\dots,S_n)$.

In what follows we will find a multivariable commutative analogue 
of this problem.
The place of the multiplication by $e^{i t}$ is taken by the   
operators $B_1,\dots,B_n$ on the symmetric  Fock space $F_s^2(H_n)$, and the
place of $C(\TT)$ is taken by  the  Toeplitz algebra
 $C^*(B_1,\dots,B_n)$.
The Fourier coefficients of a map
$\Phi$ from $C^*(B_1,\dots,B_n)$ into $B(\cH$) are given by  the evaluations
$\Phi(B^{\bf k})$, \ ${\bf k}\in \ZZ_n^+$.
 
 \begin{theorem}\label{trig-moment} Let $\Pi$ be an admissible set in  $\ZZ^+_n$ and let
 $\Gamma:\Pi \to B(\cH)$ be an  operator-valued map with $\Gamma(0)=I$. 
 Then $\Gamma$ is a 
 moment map   for a completely positive  linear map
  $\Phi:C^*(B_1,\ldots, B_n)\to B(\cH)$  if and only if
    the kernel
    $K:\Sigma_\pi\times \Sigma_\pi\to B(\cH)$ given by
  \begin{equation*}
  K(\sigma, \tau):=
  \begin{cases}
  \Gamma(\varphi(\tau\backslash \sigma)), & \text{ if } \sigma\leq \tau\\
  \Gamma(\varphi(\sigma\backslash \tau))^*, & \text{ if }\tau < \sigma\\
  0, & \text{ otherwise }
  \end{cases}
  \end{equation*}  
    is positive semidefinite.
 \end{theorem}
 \begin{proof}
 Assume $\Phi: C^*(B_1,\ldots, B_n)\to B(\cH)$  is a  completely positive 
 map  such that $\Phi(B^{\bf k})= \Gamma({\bf k})$ for any ${\bf k}\in \Pi$.
 Let  $\gamma: C^*(S_1,\ldots, S_n)\to C^*(B_1,\ldots, B_n) $ be
  the Poisson transform associated with the row contraction
   $[B_1,\ldots, B_n]$, i.e., 
   $\gamma(S_\alpha S_\beta^*)= B_\alpha B_\beta^*$, 
   \ $\alpha, \beta\in \FF_n^+$.
   Then $\Phi\circ \gamma: C^*(S_1,\ldots, S_n)\to B(\cH)$ is 
   a completely positive map such that $(\Phi\circ \gamma) (I)=I$ and 
   \begin{equation}\label{circ}
   (\Phi\circ \gamma) (S_\sigma)= \Phi(B_\sigma)= \Gamma(\varphi(\sigma)), \quad 
   \sigma\in \Sigma_\pi.
   \end{equation}
   By Stinespring's theorem \cite{S}, there is a Hilbert space 
   $\cK\supset \cH$ and a $*$-representation
    $\pi$  of $  C^*(S_1,\ldots, S_n) $  on $\cK$ such that
    
    \begin{equation}\label{Stine}
    (\Phi\circ \gamma) (A)=P_\cH \pi(A)|\cH,\quad A\in C^*(S_1,\ldots, S_n).
    \end{equation}
 Notice that the operators $V_i:=\pi (S_i)$, \ $i=1,\ldots, n$, are 
 isometries with orthogonal 
 ranges. Hence, and using relations \eqref{circ} and \eqref{Stine}, 
 we infer that 
 \begin{equation*}\begin{split}
 P_\cH V_\sigma^* V_\tau|\cH&=
 \begin{cases}
  P_\cH \pi(\tau\backslash \sigma)| \cH, & \text{ if } \sigma\leq \tau\\
  P_\cH \pi(\sigma\backslash \tau)^*| \cH, & \text{ if }\tau < \sigma\\
  0, & \text{ otherwise }
  \end{cases}\\
  &=K(\sigma, \tau).
  \end{split}
 \end{equation*}
 Therefore, if $f:\Sigma_\pi\to \cH$ is a finitely supported function, then
 \begin{equation*}\begin{split}
 \sum_{\sigma, \tau\in \Sigma_\pi} \langle K(\sigma,\tau) f(\tau), 
 f(\sigma)\rangle
 &= 
 \sum_{\sigma, \tau\in \Sigma_\pi} \langle V_\tau f(\tau), 
 V_\sigma f(\sigma)\rangle\\
 &=
 \|\sum_{ \tau\in \Sigma_\pi}  V_\tau f(\tau)\|^2\geq 0.
 \end{split}
 \end{equation*}
 
 Conversely, suppose the kernel $K$ is positive semidefinite and define
  a semidefinite form on $\cH^{\Sigma_\pi}$, the space of all finitely supported functions
   from $\Sigma_\pi$ to $\cH$, by
   \begin{equation}\label{inne}
   \langle f, \psi \rangle_K:= 
   \sum_{\sigma, \tau\in \Sigma_\pi} \langle K(\sigma,\tau) f(\tau), 
 \psi(\sigma)\rangle.
   \end{equation}
   Let $\cK$ be the Hilbert space completion of $\cH^{\Sigma_\pi}/ \cN$ where $\cN$ is the
    subspace of null vectors in the seminorm  determined
    by \eqref{inne}.
    Since $K(g_0, g_0)=I$, the Hilbert space $\cH$ imbeds isometrically in $\cK$
    by identifying $h$ with the equivalence class of the  function 
    $\sigma\mapsto \delta_{g_0}(\sigma)h$, where
     $\delta_{g_0}(\sigma)=1$ if $\sigma= g_0$,  and $0$  otherwise.
 For each $i=1,\ldots, n$, define $T_i:\cK\to \cK$ by
 \begin{equation}\label{ti}
 (T_i f)(t):= \sum_{\sigma \in  \Sigma_\pi} \delta_{g_i \sigma}(t) f(\sigma)
 \end{equation}
 for any $f\in \cH^{\Sigma_\pi}$ and $t\in \Sigma_\pi$.
 Notice that $T_i$ is a partial isometry and its range is contained
  in the closed span of functions  supported on
  \begin{equation*}
  \cR_i:=\{\sigma\in \Sigma_\pi:\ \sigma=g_i\tau  \text{ for } 
  \tau \in \Sigma_\pi\}.
  \end{equation*}
 The definition of $K$ implies that the subspaces $\cR_i$, \ $i=1,\ldots, n$, 
 are pairwise orthogonal. Therefore, 
 $[T_1,\ldots, T_n]$ is a row contraction.
 Define 
 $$
 \cE:= \bigvee_{\tau\in \FF_n^+, ~i,j=1,\ldots, n}  
 T_\tau (T_iT_j-T_jT_i)\cK
 $$
  and $\cG:= \cK\ominus \cE$.
 Since the subspace $\cE$ is invariant under each $T_i$, \ $i=1,\ldots, n$,
 it follows that $\cG$ is invariant under  each $T_i^*$, \ $i=1,\ldots, n$.
 Notice that according to \eqref{ti}, for any $f\in \cH^{\Sigma_\pi}$, we have
 \begin{equation}\label{commut}
 [(T_iT_j-T_jT_i)f](t)= \sum_{\sigma\in \Sigma_\pi}
  [\delta_{g_ig_j\sigma}(t)- \delta_{g_jg_i\sigma}(t)] f(\sigma).
 \end{equation}
 Hence, if  $f\in \cH^{\Sigma_\pi}$ and $h\in \cH$, we have
 \begin{equation*}\begin{split}
 \langle h, T_\tau (T_iT_j&-T_jT_i)f\rangle_K\\
 &= 
 \sum_{\omega, t\in \Sigma_\pi}\left< K(t, \omega) \delta_{g_0}(\omega) h, 
 \sum_{\sigma\in \Sigma_\pi}
  [\delta_{\tau g_ig_j\sigma}(t)- \delta_{\tau g_jg_i\sigma}(t)] 
  f(\sigma)\right>_\cH\\
 &= \sum_{\sigma\in \Sigma_\pi}\langle 
 K(\tau g_ig_j\sigma, g_0) h, f(\sigma)\rangle_\cH
 - \sum_{\sigma\in \Sigma_\pi}\langle 
 K(\tau g_jg_i\sigma, g_0) h, f(\sigma)\rangle_\cH,
 \end{split} 
 \end{equation*}
 if both 
 $\tau g_ig_j\sigma $
  and  $\tau g_jg_i\sigma$ are in $ \Sigma_\pi$, and $0$ otherwise.
 Due to the definition of the set  $\Sigma_\pi$, we have that
  $\tau g_ig_j\sigma\in \Sigma_\pi$
 if and only if $\tau g_jg_i\sigma\in \Sigma_\pi$. Moreover, in this case,
 due to the definition of $K$, we  get
 $K(\tau g_ig_j\sigma, g_0)=K(\tau g_jg_i\sigma, g_0)$.
 Hence, $\langle h, T_\tau (T_iT_j-T_jT_i)f\rangle_K=0$ for any 
 $\tau\in \FF_n^+$ and $f\in \cH^{\Sigma_\pi}$.
 Therefore, $h$ is orthogonal to $\cE$, which shows that $\cH\subset \cG$.
 
 Define $T_i':= P_\cG T_i|\cG$, \ $i=1,\ldots, n$. Since $\cG$ is 
  invariant under  each $T_i^*$, \ $i=1,\ldots, n$, and $[T_1,\dots, T_n]$
  is a row contraction,  we deduce that 
  $[T_1',\dots, T_n']$
  is also  a row contraction and  
  $$
  T_i'T_j'=T_j'T_i', \quad \text{ for any } i,j=1,\ldots, n.
  $$
 Now, for any $\sigma\in \Sigma_\pi$ and $h,h'\in \cH\subset\cG$, we have
 \begin{equation}\label{tprime}\begin{split}
 \langle T_\sigma'h, h'\rangle_K &= \langle  T_\sigma h, h'\rangle_K\\
 &=\sum_{\omega, t\in \Sigma_\pi}\left< K(t, \omega) \delta_{\sigma}(\omega) h,
 \delta_{g_0}(t) h' \right>_\cH\\
 &= \langle  K(g_0, \sigma) h, h'\rangle_\cH\\
 &= \langle  \Gamma(\varphi (\sigma)) h, h'\rangle_\cH.
 \end{split}
 \end{equation}
 Let $\mu:C^*(B_1,\ldots, B_n)\to B(\cG)$ be the Poisson transform associated with 
 the row contraction $[T_1',\dots, T_n']$, i.e., 
 \begin{equation}
 \label{mub}
 \mu(B^{\bf k} {B^{\bf m}}^*)
 ={T'}^{\bf k} {{T'}^{\bf m}}^*, \text{ for any } {\bf k}, {\bf m}\in \ZZ_n^+.
 \end{equation}
 Let $\Phi:C^*(B_1,\ldots, B_n)\to B(\cH)$ be the compression 
 $\Phi(A):=P_\cH \mu(A)|\cH$,\ $A\in C^*(B_1,\ldots, B_n)$. Clearly $\Phi$ is a
 unital completely positive 
 linear map.  Using relations \eqref{tprime} and \eqref{mub},
 we deduce that
 \begin{equation*}\begin{split}
 \langle \Phi(B^{\bf k})h, h'\rangle_\cH &=
  \langle P_\cH \mu(B^{\bf k})h, h'\rangle_\cH\\
  &=\langle {T'}^{\bf k}h, h'\rangle_K \\
  &= \langle  \Gamma ({\bf k}) h, h'\rangle_\cH,
  \end{split}
 \end{equation*}
 for any ${\bf k}\in \Pi$ and $h,h'\in \cH$, where $B:=(B_1,\ldots, B_n)$.
 Therefore, $ \Phi(B^{\bf k})= \Gamma ({\bf k})$ for any ${\bf k}\in \Pi$.
 The proof is complete. 
 \end{proof}
 
 Consider $\lambda_{ji}\in \CC$, $1\leq i<j\leq n$,  and let $J_\lambda$ be the 
 WOT-closed, two-sided ideal of  $F_n^\infty$ generated by 
 $\{e_j\otimes e_i-\lambda_{ji} e_i\otimes e_j:\ 1\leq i<j\leq n\}$.
 Denote $\cN_{J_\lambda}:= F^2(H_n)\ominus J_\lambda F^2(H_n)$, and
  $B_i:= P_{\cN_{J_\lambda}} S_i |\cN_{J_\lambda}$, ~$i=1,\ldots, n$.
 According to \cite{ArPo2} (see Example 3.3), if $[T_1,\ldots, T_n]$ is
  a row contraction satisfying the commutation relations 
  \begin{equation}\label{comrel}
  T_jT_i=\lambda_{ji} T_iT_j, \quad  1\leq i<j\leq n,
  \end{equation}
 then there is a unique Poisson transform $\Phi :C^*(B_1,\ldots, B_n)\to B(\cH)$
 such that $\Phi(B_\alpha B_\beta^*)=T_\alpha T_\beta^*$, \ 
 $\alpha,\beta\in \FF_n^+$. Therefore, $[B_1,\ldots, B_n]$ is the
  universal row contraction satisfying \eqref{comrel}.
 Notice that, by \eqref{comrel}, 
 for any ${\bf m}\in \ZZ_n^+$ and $\sigma\in \FF_n^+$
  such that $\varphi(\sigma)= {\bf m}$,
 there is a unique ``signature''
  $\epsilon(\sigma)\in \CC$
   (which can be expressed in terms of $\lambda_{ji}$ \cite{ArPo2}) such that 
   $T_\sigma= \epsilon(\sigma) T^{\bf m}$.
   Now define
   $$
   L_\Gamma(\sigma):= \epsilon(\sigma) \Gamma(z^{\bf m}), 
   \quad \text{where } {\bf m}\in \ZZ_n^+, ~\sigma\in \FF_n^+, ~\varphi(\sigma)= {\bf m}. 
   $$
Using this definition for $ L_\Gamma$, all the results of this section
 can be extended to  this more general setting,   with slight adjustments
 of the proofs.
 Let us remark that if $\lambda_{ji}=1$  (resp.~$\lambda_{ji}=-1$),
 \ $1\leq i<j\leq n$, then $B_1,\ldots,
 B_n$ are the creation operators on the symmetric (resp.~anti-symmetric)
  Fock space.

 \bigskip 
 
\section{Moment problems  for quotients of the free 
semigroup algebra $\CC \FF_n^+$}
 \label{extension}

Let $\cP_h$ be a set of homogeneous polynomials in  $ \CC\FF_n^+$. 
We recall that  
 $ p=\sum_{j=1}^m a_j {\alpha_j}$, $a_j\neq 0$,  is a homogeneous polynomial if 
  $|\alpha_j|=|\alpha_1|$ for any $j=1,\ldots, m$.
Let $J$ be the WOT-closed, two-sided ideal of $F_n^\infty$ generated by
the polynomials  $\{p(S_1,\ldots, S_n):\ p\in \cP_h\}$.
Denote $\cN_{J}:= F^2(H_n)\ominus J F^2(H_n)$, and define
  $B_i:= P_{\cN_{J}} S_i |\cN_{J}$, ~$i=1,\ldots, n$. 
  Throughout this section,  $C^*(B_1,\ldots, B_n)$ (resp. $\cA_n$)
  is the $C^*$-algebra (resp. non-selfadjoint norm closed algebra)
   generated by $B_1,\ldots, B_n$ and  the identity. 

Using the results from \cite{Po-poisson} and \cite{ArPo2},  we can prove 
the following result. We only sketch the proof.
\begin{theorem}\label{poisson}
A map $\mu: C^*(B_1,\ldots, B_n)\to B(\cH)$ is a Poisson transform 
if and only if there is a  row contraction $[T_1,\ldots, T_n]$, 
$T_i\in B(\cH)$, with
$p(T_1,\ldots, T_n)=0$ for any homogeneous polynomial $p\in \cP_h$, such that
  $$
  \mu(B_\alpha B_\beta^*)=T_\alpha T_\beta^*, \quad \alpha,\beta\in \FF_n^+.
  $$
\end{theorem}
\begin{proof}
Suppose $[T_1,\ldots, T_n]$, 
$T_i\in B(\cH)$,  is a row contraction with
$p(T_1,\ldots, T_n)=0$ for any  $p\in \cP_h$.
Let $K_r(T)$ be the Poisson kernel associated with $T:=[T_1,\ldots, T_n]$, 
as defined in the introduction.
Then, for any homogeneous polynomial
$p= \sum_{j=1}^m a_j \alpha_j \in \cP_h$, ~$\omega, 
\beta\in \FF_n^+$, and $k,h\in \cH$, we have
$$
\langle K_r(T) k, S_\omega p(S_1,\ldots, S_n) e_\beta\otimes h\rangle
=r^{|\alpha_1|+|\omega|+|\beta|}
\langle  k, T_\omega p(T_1,\ldots, T_n) T_\beta h\rangle=0.
$$
Since  the noncommutative  analytic Toeplitz algebra $F_n^\infty$ is the WOT-closure of the polynomials in $S_1,\ldots, S_n, I$, it is easy to see
that
 $$
\langle K_r(T) k, fp(S_1,\ldots, S_n) g\otimes h\rangle=0
$$
for any $f\in F_n^\infty, g\in F^2(H_n)$, and $p\in \cP_h$. Now, it is clear that 
$K_r(T) k\in \cN_J\otimes \cH$ for any $k\in \cH$ and $0<r<1$. 
According to Remark 3.2 from
 \cite{ArPo2}, there is a Poisson transform $\mu$ satisfying the required conditions. 
The converse is straightforward.
\end{proof}

We say that a set $\Sigma\subset \FF_n^+$ is compatible with  
$p=\sum_{j=1}^m a_j\alpha_j\in \CC\FF_n^+$ if, for any  
 $\omega, \beta\in \FF_n^+$, we have either
 \begin{equation}\label{compat}
 \{\omega \alpha_1 \beta,\ldots, \omega \alpha_n \beta\}\subset \Sigma,\ 
\text{  or }\ 
\{\omega \alpha_1 \beta,\ldots, \omega \alpha_n \beta\}\cap \Sigma=\emptyset.
\end{equation} 
A pair $(\Sigma, \Lambda)$, $\Lambda\subset \CC\FF_n^+$,
 is called admissible if $\Sigma$ is an admissible  set in 
$\FF_n^+$ and compatible with each  $p\in \Lambda$.
\begin{remark} 
 If $\Sigma=\FF_n^+$ or $\Sigma=\{\alpha\in \FF_n^+:\ |\alpha|\leq k\}$, 
 \ $k\in \ZZ^+$, then  the pair $(\Sigma, \cP_h)$ is admissible for any 
 set $\cP_h\subset \CC\FF_n^+$ of homogeneous polynomials. Moreover,
 $(\FF_n^+, \Lambda)$ is admissible for any  set $\Lambda\subset \CC\FF_n^+$.
\end{remark}

Let $\cP_h\subset \CC\FF_n^+$ be a set of homogeneous polynomials and 
let $J_h$ be the two-sided ideal of $\CC\FF_n^+$ generated by $\cP_h$, i.e.,
$$
J_h=\{\sum_{j=1}^k x_i p_i y_i:\ x_i, y_i\in \FF_n^+, ~p_i\in \cP_h, ~k\geq 1\}.
$$
Consider the quotient algebra $\CC\FF_n^+/ J_h$ and denote by $\hat{p}$ 
the coset corresponding to 
$p\in \CC\FF_n^+$, i.e., $\hat{p}= p+J_h$. Notice that if $q$ is 
any representative of the coset $\hat{p}$, then
 $q(B_1,\ldots, B_n)=p(B_1,\ldots, B_n)$, where $B_1,\ldots, B_n$ are 
 defined as above. Hence, the map $\Phi: \CC\FF_n^+/ J_h\to \cA_n$ defined by
 $\Phi(\hat{p}):= p(B_1,\ldots, B_n)$ is a homomorphism.
\begin{theorem}\label{main-general}
Let  $(\Sigma, \cP_h)$ be an admissible pair and let $\hat{L}: 
\Sigma/J_h\to B(\cH)$
be an operator-valued map.  Then there is a Poisson transform 
$\rho: C^*(B_1,\ldots, B_n)\to B(\cH)$ such that
\begin{equation}\label{hat}
\hat{L}(\hat{p})=\rho(p(B_1,\ldots, B_n)) \hat{L}(\hat{g_0}), \quad \hat{p}\in 
\CC \Sigma/J_h,
\end{equation}
if and only if the following conditions hold:
\begin{enumerate}
\item[(i)] If $p=\sum_{j=1}^m a_j\alpha_j\in \cP_h$ and 
 $\omega, \beta\in \FF_n^+$ are such that 
 $\{\omega \alpha_1 \beta,\ldots, \omega \alpha_m \beta\}\subset \Sigma$, then
\begin{equation}\label{eq}
\sum_{j=1}^m a_j\hat{L}(\widehat{\omega\alpha_j \beta})=0.
\end{equation}
\item[(ii)]
$K_1\leq K_2$,  where the kernels
 $K_1$ and $K_2$ are associated, as in  \eqref{k1} and \eqref{k2},
  with the admissible  set
  $\Sigma $ and the map $L:\Sigma\to B(\cH) $ defined  by 
  $L(\sigma):= \hat{L}(\hat{\sigma})$, 
  \ $\sigma\in \Sigma$.

\end{enumerate}
\end{theorem}
  \begin{proof} 
  Let $\rho: C^*(B_1,\ldots, B_n)\to B(\cH)$ be a Poisson transform 
  satisfying \eqref{hat}. Then, for any $(\alpha, \beta)\in \Lambda_\Sigma$,
  we have  
  $$
  L(\alpha \beta)=\hat{L}(\widehat{\alpha \beta})= \rho(B_\alpha )
  \hat{L}(\hat{\beta})= \mu(S_\alpha)L(\beta),
  $$
  where $\mu:C^*(S_1,\ldots, S_n)\to B(\cH)$ is the Poisson transform associated to
  the row contraction
  $[\rho(B_1),\ldots, \rho(B_n)]$. Hence,  $\mu$ is a representing Poisson transform 
  for the map $L$ and, according to Theorem \ref{main}, we have $K_1\leq K_2$.  
  On the other hand, 
  if $p=\sum_{j=1}^m a_j\alpha_j\in \cP_h$ and 
 $\omega, \beta\in \FF_n^+$ are such that 
 $\{\omega \alpha_1 \beta,\ldots, \omega \alpha_m \beta\}\subset \Sigma$,
 then $\sum_{j=1}^m a_j (\omega \alpha_j \beta)$ is in $\CC \Sigma $ and 
 relation \eqref{hat} implies
 \begin{equation*}
 \sum_{j=1}^m a_j\hat{L}(\widehat{\omega\alpha_j \beta})=
 \rho(B_\omega) \rho(\sum_{j=1}^m a_j  B_{\alpha_j}) \rho(B_\beta) \hat L(\hat g_0).
 \end{equation*}
 Since $\sum_{j=1}^m a_j  B_{\alpha_j}=0$, relation \eqref{eq} is satisfied.
 
 Conversely, as in the proof of Theorem \ref{main}, we  can find  the operators
 $T_i\in B(\cH)$, \ $i=1,\ldots, n$, with the following properties:
 \begin{enumerate}
 \item[(i)] $[T_1,\ldots, T_n]$ is a contraction;
 \item[(ii)] 
  $
 T_iL(\sigma)h=\begin{cases}L(g_i\sigma)h, & \text{ if } g_i\sigma\in \Sigma\\
 0, & \text{ otherwise};
 \end{cases}
 $
 \item[(iii)] $T_i|\cM^\perp=0$, where $\cM$ is the closed linear span of all
  vectors $L(\sigma)h$, where $\sigma\in \Sigma$ and $h\in \cH$.
 \end{enumerate}
 Let 
  $p=\sum_{j=1}^m a_j\alpha_j\in \cP_h$, $\sigma\in \Sigma$, and $h\in \cH$.
 Since $(\Sigma, \cP_h)$ is an admissible pair,  we have  either 
 $
 \{\alpha_1\sigma ,\ldots,  \alpha_m \sigma\}\subset \Sigma$,\ 
 or 
$\{ \alpha_1\sigma ,\ldots,  \alpha_m \sigma\}\cap \Sigma=\emptyset$.
Using the above properties,  we obtain
$$
 \sum_{j=1}^m a_jT_{\alpha_j} L(\sigma)h=
 \begin{cases}  \sum\limits_{j=1}^m a_j L(\alpha_j\sigma)h, & \text{ if }
  \{\alpha_1\sigma ,\ldots,  \alpha_m \sigma\}\subset \Sigma\\
 0, & \text{ otherwise};
 \end{cases}
 $$
 Using relation \eqref{eq} and the property (iii), we deduce 
 $\sum_{j=1}^m a_jT_{\alpha_j}=0$, i.e., $p(T_1,\ldots, T_n)=0$ for 
 any $p\in \cP_h$.
 According to Theorem \ref{poisson}, there is a Poisson transform
 $\rho: C^*(B_1,\ldots, B_n)\to B(\cH)$ such that
 $\rho(B_\alpha B_\beta^*)= T_\alpha T_\beta^*$, \ $\alpha,\beta\in \FF_n^+$.
 Now, it is easy to see that, if $\omega\in \Sigma$, then
 $$
 \hat{L}(\hat{\omega})=L(\omega) =T_\omega L(g_0)= \rho(B_\omega)\hat{L}(\hat{g_0}).
 $$
 Therefore, relation \eqref{hat} is satisfied. This completes the proof.
\end{proof}
As in Section \ref{Moment}, one can obtain a vector-valued version of Theorem
\ref{main-general}. This provides a characterization 
   for  the orbits of contractive  
Hilbert modules over  the quotient algebra $\CC \FF_n^+/J_h$,
 where $J_h$ is  the 
two-sided ideal  generated by $\cP_h$. A more general result is presented 
in Section
\ref{moment-Cuntz}.


 
 In what follows, 
  we solve the operator-valued trigonometric moment problem for 
completely positive linear maps on the $C^*$-algbras 
    $C^*(B_1,\ldots, B_n)$ described above.
    
 \begin{theorem}\label{trig2} Let  $(\Sigma, \cP_h)$ be an admissible pair and 
 let $\hat{L}: \Sigma/J_h\to B(\cH)$
be an operator-valued map with 
 $\hat{L}(\hat{g_0})=I$.  
Then there is a   completely positive linear map 
$\Phi: C^*(B_1,\ldots, B_n)\to B(\cH)$ such that
\begin{equation}\label{hat1}
\hat{L}(\hat{p})=\Phi(p(B_1,\ldots, B_n)), \quad \hat{p}\in 
\CC \Sigma/J_h,
\end{equation}
if and only if the following conditions hold:
\begin{enumerate}
\item[(i)] If $\sum_{j=1}^m a_j\alpha_j\in \cP_h$ and 
 $\omega, \beta\in \FF_n^+$ are such that 
 $\{\omega \alpha_1 \beta,\ldots, \omega \alpha_m \beta\}\subset \Sigma$, then
$$
\sum_{j=1}^m a_j\hat{L}(\widehat{\omega\alpha_j \beta})=0.
$$
\item[(ii)]
The kernel
    $K:\Sigma\times \Sigma\to B(\cH)$ given by
  \begin{equation*}
  K(\sigma, \tau):=
  \begin{cases}
  L(\tau\backslash \sigma), & \text{ if } \sigma\leq \tau\\
  L(\sigma\backslash \tau)^*, & \text{ if }\tau < \sigma\\
  0, & \text{ otherwise }
  \end{cases}
  \end{equation*}  
    is positive semidefinite, where the  map $L:\Sigma\to B(\cH) $
      is defined  by 
  $L(\sigma):= \hat{L}(\hat{\sigma})$, 
  \ $\sigma\in \Sigma$.
\end{enumerate}
 \end{theorem}
 \begin{proof}
 Assume there is a   completely positive linear map 
$\Phi: C^*(B_1,\ldots, B_n)\to B(\cH)$ such that \eqref{hat1} holds.
If $p=\sum_{j=1}^m a_j\alpha_j\in \cP_h$, then $p(B_1,\ldots, B_n)=0$.  
 Moreover, if $\omega, \beta\in \FF_n^+$ 
 are such that 
 $\{\omega \alpha_1 \beta,\ldots, \omega \alpha_m \beta\}\subset \Sigma$, then
 $$
 \sum_{j=1}^m a_j\hat{L}(\widehat{\omega\alpha_j \beta})=
 \Phi[B_\omega p(B_1,\ldots, B_n) B_\beta]=0.
 $$
 On the other hand, let  
 $\gamma: C^*(S_1,\ldots, S_n)\to C^*(B_1,\ldots, B_n) $ be
  the Poisson transform associated with the row contraction
   $[B_1,\ldots, B_n]$.
   Then $\Phi\circ \gamma: C^*(S_1,\ldots, S_n)\to B(\cH)$ is 
   a completely positive map such that $(\Phi\circ \gamma) (I)=I$ and 
    $$
   (\Phi\circ \gamma) (S_\sigma)=  L(\sigma), \quad 
   \sigma\in \Sigma.
    $$
 Now, as in the proof of Theorem \ref{trig-moment}, we can prove that the kernel $K$
 is positive semidefinite.
 
 Conversely, assume $K\geq 0$ and relation \eqref{eq} holds. This part of the proof
 is similar to that of Theorem \ref{trig-moment}. Let us follow that proof and mention 
 only the  necessary changes.
 First, the subspace $\cE$ should be defined by $$
 \cE:=\bigvee_{\tau\in \FF_n^+, ~p\in \cP_h} T_\tau  p(T_1,\ldots, T_n)\cK.
 $$
 According to \eqref{ti}, for any $f\in \cH^\Sigma$ and
$p=\sum_{j=1}^m a_j\alpha_j\in \cP_h$, we have
$$
[p(T_1,\ldots, T_n) f](t)= \sum_{\sigma\in \Sigma} 
\left[\sum_{j=1}^m a_j \delta_{\alpha_j \sigma}(t)\right] f(\sigma).
$$  
 Hence, if $\tau\in \FF_n^+$ and $h\in \cH$, then we deduce
 \begin{equation*}
 \begin{split}
 \langle h, T_\tau  p(T_1,\ldots, T_n)f\rangle_\cK &=
 \sum_{\sigma\in \Sigma} 
\sum_{j=1}^m \bar{a_j} \langle K(\tau \alpha_j \sigma, g_0)h, f(\sigma)\rangle_\cH\\
&= 
\sum_{\sigma\in \Sigma} \left<
\sum_{j=1}^m \bar{a_j}  L(\tau \alpha_j \sigma)^*h, f(\sigma)\right>_\cH,
 \end{split}
 \end{equation*}
 if $\{\omega \alpha_1 \beta,\ldots, \omega \alpha_m \beta\}\subset \Sigma$, and 
 $0$ otherwise. Here we used the fact that the pair $(\Sigma, \cP_h)$ is admissible.
 By \eqref{hat1}, we have $\sum_{j=1}^m a_j L(\tau \alpha_j \sigma)=0$. Therefore,
 $\cH\subset \cG:=\cK\ominus \cE$.
 Now, define $T_i':= P_\cG T_i |\cG$, \ $i=1,\ldots, n$, and notice that
 $[T_1',\ldots, T_n']$ is a row contraction such that  $p(T_1',\ldots, T_n')=0$ 
 for any $p\in \cP_h$.
 According to Theorem \ref{poisson}, there is a Poisson transform
 $\mu: C^*(B_1,\ldots, B_n)\to B(\cG)$ associated with $[T_1'\ldots, T_n']$.
 Let $\Phi:C^*(B_1,\ldots, B_n)\to B(\cH)$ be the compression 
 $\Phi(A):=P_\cH \mu(A)|\cH$,\ $A\in C^*(B_1,\ldots, B_n)$. It is clear 
 that $\Phi$ is a
 unital completely positive 
 linear map.   
   On the other hand,  it is easy to see that
 \begin{equation*}\begin{split}
 \langle \Phi(B_\alpha)h, h'\rangle_\cH &=
  \langle P_\cH \mu(B_\alpha)h, h'\rangle_\cH =\langle T_\alpha h, h'\rangle_K \\
  &=\sum_{\omega, t\in \Sigma_\pi}\left< K(t, \omega) \delta_{\alpha}(\omega) h,
 \delta_{g_0}(t) h' \right>_\cH\\
  &=\langle K(g_0, \alpha) h, h' \rangle
  = \langle  L (\alpha) h, h'\rangle_\cH\\
   &=\langle  \hat{L} (\hat{\alpha}) h, h'\rangle_\cH
  \end{split}
 \end{equation*}
 for any $\alpha\in \Sigma$ and $h,h'\in \cH$. 
 The proof is complete. 
 \end{proof}

 \bigskip

\section{Moment problems for the free product semigroup $*_{i=1}^n G_i^+$}
 \label{moment-Cuntz}


Let $G_i$ be any discrete subgroups of the real line $\RR$, and let $G_i^+$
denote the positive cones $G_i^+:= G_i \cap [0,\infty)$. We assume that at least
one subgroup $G_i$ is dense in $\RR$. 
The left
regular representation of the free product semigroup   $*_{i=1}^n G_i^+$
   is a
semigroup of isometries $\lambda(g)$ on the Hilbert space 
$\ell^2(*_{i=1}^n G_i^+)$
 with orthonormal
basis $\{\xi_g:g\in  *_{i=1}^n G_i^+\}$, given by 
\[  \lambda(g)\xi_\sigma = \xi_{g\sigma}, \quad g,\sigma \in *_{i=1}^n G_i^+.\] 
 The reduced semigroup $C^*$-algebra $C_r^*(*_{i=1}^n G_i^+)$  is 
the C*-algebra generated
by the left regular representation of $*_{i=1}^n G_i^+$. 
In this new setting, the noncommutative disc algebra 
 $\cA:=\cA(*_{i=1}^n G_i^+)$ 
   is defined
  \cite{DPo} as the
norm-closed algebra generated by the left regular representation.

   Consider
representations $T_i$ of $G_i^+$ as semigroups of contractions on a common
Hilbert space $\cH$ such that
\begin{equation}\label{contrep}
  \sum_{i=1}^n T_i(g_i) T_i(g_i)^* < I,  \qquad g_i \in G_i^+\backslash\{0\} .
\end{equation}
Such  representations (uniquely) determine a contractive representation of the
free product semigroup $*_{i=1}^n G_i^+$ by sending an element
$\sigma=g_1g_2\dots g_k$, where $g_j\in G_{i_j}^+$, to the  contraction
\[
  T(\sigma) := T_{i_1}(g_1)T_{i_2}(g_2)\dots T_{i_k}(g_k),
\]
and $T(g_0)=I$, where $g_0$ denotes the neutral element
 of $*_{i=1}^n G_i^+$.
 The generalized Cuntz algebra $\cO(*_{i=1}^n G_i^+)$ was defined in \cite{DPo}
  as the universal $C^*$-algebra generated by isometric representations $V_i$ of 
  the semigroups $G_i^+$  satisfying condition \eqref{contrep}. Extending some 
  results of Douglas \cite{Do}  and Cuntz \cite{Cu}, it was proved in \cite{DPo}
  (see also  \cite{D1}, 
   \cite{D2}, \cite{Ni},
    \cite{DPo}, \cite{LR}, \cite{L}) that if $n\geq 2$, then the $C^*$-algebra
    generated by $V_i$ does not depend
    on the choice of the isometric  representations $V_i$. Therefore
  $\cO(*_{i=1}^n G_i^+)=C_r^*(*_{i=1}^n G_i^+)$. When $n=1$, we set
   $\cO(G^+):= C_r^*(G^+)$.

A map $\mu:\cO(*_{i=1}^n G_i^+)\to B(\cH)$ is called Poisson transform on the
generalized Cuntz algebra $\cO(*_{i=1}^n G_i^+)$ on the Hilbert space $\cH$ 
if it is a unital completely positive
 linear  map such that 
    $$
    \mu(AX)=\mu(A)\mu(X), \quad \text{ for any } X\in  \cA\cA^*
    \text{ and }
A\in \cA.
$$
  According to \cite{DPo}, a map $\mu$ is a Poisson transform on
  $\cO(*_{i=1}^n G_i^+)$ if and only if 
  there are 
   contractive representations
   $T_i$  
    of $G_i^+$ on
a common Hilbert space $\cH$ satisfying the  condition \eqref{contrep}, 
and  
  such that 
\[
  \mu(\lambda(\sigma) \lambda(\tau)^*) = T_\sigma T_\tau^*, \quad 
  \sigma,\tau \in *_{i=1}^n G_i^+.
\]
As in the previous sections, a set $\Sigma\subset  *_{i=1}^n G_i^+$ is 
called admissible if
$\alpha \beta\in \Sigma$ implies $\beta\in \Sigma$.
The following result is an extension of Theorem \ref{main} to generalized Cuntz
algebras.

  \begin{theorem}\label{main-semigroup}
  Let $\Sigma\subseteq *_{i=1}^n G_i^+ $ be an admissible
   set and let $L:\Sigma\to B(\cH)$ be
   an operator-valued  map.
   Then there is a Poisson transform  
  $\mu:\cO(*_{i=1}^n G_i^+)\to B(\cH)$  
  such that
  \begin{equation}\label{L}
  L(\sum a_\sigma \sigma)=\mu(\sum a_\sigma \lambda(\sigma))L(g_0)
  \end{equation}
  for any  polynomial \ $\sum a_\sigma \sigma \in \CC\Sigma$
   if and only if $K_1\leq K_2$, where the kernels $K_1, K_2$ are defined
   as  in \eqref{k1} and \eqref{k2}.
  \end{theorem}
\begin{proof}
Assume that $\mu:\cO(*_{i=1}^n G_i^+)\to B(\cH)$  is a Poisson
 transform satisfying \eqref{L} and
  let $T_i(g_i):= \mu(\lambda(g_i)), g_i\in G_i^+$.
Notice that $T_i$ are contractive representations of $G_i^+$ on $\cH$ satisfying
 condition \eqref{contrep}.
According to \cite{DPo} (see also \cite{SzF-book}, \cite{Ml},  and 
\cite{Po-isometric}
for important particular cases obtained before), there is a Hilbert space $\cK$ containing $\cH$
 and isometric
representations $V_i$ of $G_i^+$ on $\cK$ such that
\begin{enumerate}
\item[(i)] $\sum_{i=1}^n V_i(g_i) V_i(g_i)^* < I$ 
  for all $g_i \in G_i^+\backslash\{0\}$.
\item[(ii)] $V_i(g_i)^*|_\cH = T_i(g_i)^*$
  for all $g_i \in G_i^+$, $1\le i\le n$.
\end{enumerate}
Notice that if $V$ is the representation of $*_{i=1}^n G_i^+$ associated
 with $V_i$,
and $\sigma,\tau\in *_{i=1}^n G_i^+$ are two elements such that neither
 is a multiple of the other, then $V(\sigma)$ and $V(\tau)$ are isometries 
 with orthogonal ranges. 
 Now, the proof of the direct implication  is an extension of that of  
    Theorem \ref{main}.

   For the converse, we only sketch
  a proof.  Assume $K_1\leq K_2$. The representations $V_i$ of $G_i^+$ 
  are defined by
   $V_i(g_i)f=\varphi$, where $\varphi$ is given by \eqref{var}. 
It is easy to see that $V_i(g_i)$ is a partial isometry and
 $V_i(g_i)V_i(g_i')=V_i(g_i')V_i(g_i)$
 for any $g_i, g_i'\in G_i^+$.
 Moreover, $V_i( G_i\backslash \{0\})$  has range contained 
 in the closed span of functions supported
 on
 $$
 \cR_i:=\{(\alpha,\beta)\in \Lambda_\Sigma:\ \alpha=g_i\tau \text{ for } 
 g_i\in G_i\backslash \{0\}, \tau\in *_{i=1}^n G_i^+\}.
 $$
 The definition of $K_2$ guarantees that the  
 subspaces $\cR_i$ are pairwise orthogonal. Hence 
 the representations $V_i$, $i=1,\ldots, n$,  satisfy condition \eqref{contrep}.
 Define 
  the representation $ T_i$ of $G_i^+$  on $\cH$ by
  $$
  T_i(g_i):= XV_i(g_i)X^*\quad \text{ if }
   g_i\in G_i^+\backslash \{0\}, 
  $$
   and $T_i(0):=I$,
  where the operator $X$ is given by \eqref{xf}.
  Clearly the representations $ T_i$  also satisfy condition \eqref{contrep} and 
  $$T_i(g_i) L(\sigma)h= L(g_i\sigma)h
  $$
  if $g_i \sigma\in \Sigma$, $g_i\in G_i^+$, and $h\in \cH$.
  Hence $T(\omega)L(g_0)=L(\omega)$, \ $\omega\in \Sigma$, where $T$ is the
  contractive  representation 
  of $*_{i=1}^n G_i^+$ associated with $T_i$.
 Let  $\mu:\cO(*_{i=1}^n G_i^+)\to B(\cH)$  be the  Poisson
 transform associated with $T_i$. Since 
 $\mu(\lambda(\sigma))= 
 T(\sigma)$, \ $\sigma\in *_{i=1}^n G_i^+$,  we infer that 
 $L(\omega)=\mu(\lambda(\omega)) L(g_0)$, 
 \ $\omega\in \Sigma$, and the proof is complete.

\end{proof}

We remark that all the results of Section \ref{Moment}
hold true  in this new setting. To get the analogues just replace 
$C^*(S_1,\ldots, S_n)$
by the generalized Cuntz algebra $\cO(*_{i=1}^n G_i^+)$. 
The proofs remain essentially the same. Let us mention that an analogue
of the classical trigonometric moment problem for
$\cO(*_{i=1}^n G_i^+)$ was obtained in \cite{DPo}.


%
Concerning the commutative case, 
let us first remark that there are nontrivial contractive representations
 $T:\prod_{i=1}^n G_i^+\to B(\cH)$  such that
 $T_i:=T| G_i^+$ satisfies condition \eqref{contrep}.
Indeed, define the subspace $\cM$ of $\ell^2(*_{i=1}^n G_i^+)$ by setting
$$
\cM:= \bigvee_{\omega,\sigma,\tau\in *_{i=1}^n G_i^+}
\lambda(\omega)[\lambda(\sigma) \lambda(\tau)-\lambda(\tau)\lambda(\sigma)]
\ell^2(*_{i=1}^n G_i^+),
$$
where ~$\lambda$ is 
 the left
regular representation of the free product semigroup   $*_{i=1}^n G_i^+$
    on $\ell^2(*_{i=1}^n G_i^+)$.
Notice that the subspace  $\cN:= \ell^2(*_{i=1}^n G_i^+)\ominus \cM$ is
  invariant under each 
$\lambda(\omega)^*$,  ~$\omega\in *_{i=1}^n G_i^+$.
Let  $\Psi:\prod_{i=1}^n G_i^+\to B(\cN)$  be given by 
$$
\Psi({\bf g}):= P_\cN \lambda(g_1)\cdots \lambda(g_n)|\cN,\quad 
{\bf g}=(g_1,\ldots, g_n)\in \prod_{i=1}^n G_i^+.
$$
It is clear now that $\Psi$  is a nontrivial representation satisfying
the inequality \eqref{contrep}.
Notice that in the particular case when
$G_i^+=\NN$, \ $i=1,\ldots, n$, the $C^*$-algebra generated by
$
\{\Psi({\bf g}):\ {\bf g}\in \prod_{i=1}^n G_i^+\}
$ coincides with the Toeplitz algebra $C^*(B_1,\ldots, B_n)$,  where
$B_1,\ldots, B_n$ are the creation operators on the symmetric Fock space.

  The commutative semigroup $\prod_{i=1}^n G_i^+ $ is endowed with 
  the product order $\ll$. We denote by ${\bf g}_0$ its neutral element.
 A   set  $\Pi\subset \prod_{i=1}^n G_i^+ $  is called 
   admissible
  if ${\bf g}\in \Pi$ and
  ${\bf f} \ll {\bf g}$ imply ${\bf f}\in \Pi$.
  There is a canonical homomorphism
$\varphi$ of $*_{i=1}^n G_i^+ $ onto $\prod_{i=1}^n G_i^+ $ such that 
it is the identity on  each $G_i^+$.
Notice that if $\Pi$ is an admissible set for $\prod_{i=1}^n G_i^+$, then
\begin{equation*}
\Sigma_\pi:= \{\sigma\in *_{i=1}^n G_i^+ :\ \varphi(\sigma)\in \Pi\}
\end{equation*}
is an admissible set for $*_{i=1}^n G_i^+ $. 
We associate with each map $\Gamma:\Pi \to B(\cH)$ 
another  map $L_\Gamma:\Sigma_\pi\to B(\cH)$ by setting $L_\Gamma(\sigma):= 
\Gamma (\varphi(\sigma))$, 
$\sigma\in \Sigma_\pi$.

  \begin{theorem}\label{commut-semigroup}
  Let $\Pi\subseteq \prod_{i=1}^n G_i^+ $ be an admissible
   set and let $\Gamma:\Pi\to B(\cH)$.
   Then there is  a contractive  representation $T:\prod_{i=1}^n G_i^+ \to B(\cH)$
   satisfying  the   condition \eqref{contrep} such that
  \begin{equation}\label{ga}
  \Gamma({\bf g})=T({\bf g}) \Gamma ({\bf g}_0), \quad  {\bf g}\in \Pi,
  \end{equation}
    if and only if 
     $K_1\leq K_2$, 
     where the kernels
 $K_1$ and $K_2$ are associated with the admissible  set
  $\Sigma_\pi\subset *_{i=1}^n G_i^+ $ and the map $L_\Gamma:\Sigma_\pi\to B(\cH)$, 
  as defined by \eqref{k1} and \eqref{k2}.
  \end{theorem}
  \begin{proof}
   Assume $T:\prod_{i=1}^n G_i^+ \to B(\cH)$ is a contractive representation 
   satisfying     conditions \eqref{contrep} and \eqref{ga}.
   As in the proof of Theorem \ref{main-semigroup}, there is an isometric 
   representation $V$ of $*_{i=1}^n G_i^+$ on a Hilbert space 
   $\cK\supset \cH$ which dilates the contractive representation 
   of $*_{i=1}^n G_i^+$ associated with
   $T_i:=T|G_i^+$, \ $i=1\ldots, n$.
   For any $(\alpha, \beta)\in \Lambda_{\Sigma_\pi}$, we have
   \begin{equation}
   \begin{split}
   L_\Gamma(\alpha\beta)&= \Gamma(\varphi(\alpha)+\varphi(\beta))\\
   &= T(\varphi(\alpha)) \Gamma(\varphi(\beta))\\
   &=P_\cH V(\alpha) L_\Gamma(\beta).
   \end{split}
   \end{equation}
Now, as in the proof of Theorem \ref{main}, one can show  that $K_1\leq K_2$.

Conversely, assume $K_1\leq K_2$. Following the proof of Theorem \ref{main}
(see also Theorem \ref{main-semigroup}), we find representations 
$T_i$ of $G_i^+$ on $\cH$ defined by
$T_i(g_i):= XV_i(g_i)X^*$, $ g_i\in G_i^+$, such that the condition 
\eqref{contrep} holds and
\begin{equation*}
 T_i(g_i)L_\Gamma(\sigma)h=
 \begin{cases}   L_\Gamma(g_i\sigma)h,
 & \text{ if }  g_i\sigma\in \Sigma_\pi\\
    0, & \text{ otherwise}.
\end{cases}
 \end{equation*}
Moreover, due to the definitions of $\Sigma_\pi, L_\Gamma$ and $K_2$,
 we infer that
 $$T_i(g_i)T_j(g_j)=T_j(g_j)T_i(g_i), \quad g_i\in G_i^+, g_j\in G_j^+.
 $$ 
Let $T:\prod_{i=1}^n G_i^+ \to B(\cH)$  be the  representation 
associated with  $T_i$.
If ${\bf g}=(g_1,\ldots, g_n)\in \Pi$ and $h\in \cH$, then we have
\begin{equation*}
\begin{split}
\Gamma({\bf g})h&= L_\Gamma(g_1\cdots g_n)h\\
&= T_1(g_1)\cdots T_n(g_n) L_\Gamma( g_0)h\\
&=T({\bf g }) \Gamma({\bf g}_0)h.
\end{split}
\end{equation*}
The proof is complete.

\end{proof}

If $\Pi\subseteq \prod_{i=1}^n G_i^+ $ is an admissible
   set,  
 we associate with any map $N:\Pi \to \cH$ 
another map $M_N:\Sigma_\pi\to \cH$ by setting $M_N(\sigma):= 
N (\varphi(\sigma))$, 
$\sigma\in \Sigma_\pi$.
As in Section \ref{commutative}, one can use 
Theorem \ref{commut-semigroup} to obtain  the following vector-valued 
moment problem.
 \begin{theorem} Let $\Pi$ be an admissible set in  $\prod_{i=1}^n G_i^+$ 
 and let
 $N:\Pi \to \cH$ be a map. 
   Then there is  a contractive representation $T:\prod_{i=1}^n G_i^+ \to B(\cH)$
   satisfying     condition \eqref{contrep} such that
  \begin{equation*}
  N({\bf g})=T({\bf g}) N ({\bf g}_0), \quad  {\bf g}\in \Pi,
  \end{equation*}
  if and only if $K_3\leq K_4$, where the kernels
 $K_3$ and $K_4$ are associated with the admissible  set
  $\Sigma_\pi\subset *_{i=1}^n G_i^+ $ and the map $M_N:\Sigma_\pi\to \cH$, 
  as defined by \eqref{k3} and \eqref{k4}.
 \end{theorem}
 %

 In what follows we solve  the operator-valued moment problem
  for Poisson transforms on the generalized Cuntz algebra
 $\cO(*_{i=1}^n G_i^+)$, satisfying polynomial identities.
 The result    provides, 
 in particular, a characterization for the orbits  of contractive Hilbert modules 
 over  the quotient algebra $\CC *_{i=1}^n G_i^+/J$, where $J$ is a  
 two-sided ideal of the free semigroup algebra $\CC *_{i=1}^n G_i^+$.
 First, let us show that, in general,  there are  nontrivial representations $\Psi$
 of $\CC *_{i=1}^n G_i^+/J$ such that
  \begin{equation}\label{psihat}
  \sum_{i=1}^n \Psi(\hat {g_i})\Psi(\hat {g_i})^*\leq I, \quad 
   g_i\in G_i^+\backslash  \{0\}.
  \end{equation}
 Define the subspace $\cM_J$ of $\ell^2(*_{i=1}^n G_i^+)$ by setting
$$
\cM_J:= \bigvee_{\omega \in *_{i=1}^n G_i^+, ~ \sum a_\sigma {\sigma}\in J}
\lambda(\omega)\left[ \sum a_\sigma \lambda({\sigma})\right]
\ell^2(*_{i=1}^n G_i^+).
$$
Clearly, the subspace  $\cN_J:= \ell^2(*_{i=1}^n G_i^+)\ominus \cM_J$ is
  invariant under each 
$\lambda(\omega)^*$,  ~$\omega\in *_{i=1}^n G_i^+$.
 Notice that if $\cN_J\neq \CC$, then 
 the map $\Psi:\CC *_{i=1}^n G_i^+/J\to B(\cN_{ J})$ given  by 
 $$
 \Psi(\widehat{\sum a_\sigma \sigma }) :=P_{\cN_{J}}
  \left(\sum a_\sigma \lambda(\sigma)\right)|\cN_{ J}
  $$ 
 is a  nontrivial representation of   
 $\CC *_{i=1}^n G_i^+/J$, satisfying \eqref{psihat}.

\begin{theorem}\label{exten}
Let  $(\Sigma, J)$ be an admissible pair 
where $\Sigma \subset *_{i=1}^n G_i^+$ and $J$  is a two-sided 
ideal of $\CC*_{i=1}^n G_i^+$.
Let $\hat{L}: \Sigma/J\to B(\cH)$
be an operator-valued map.  Then there is a Poisson transform 
$\mu: \cO(*_{i=1}^n G_i^+)\to B(\cH)$ such that
\begin{enumerate}
\item[(i)] $\mu(\sum a_\sigma \lambda(\sigma))=0$ for any 
~$ \sum a_\sigma \sigma\in J$, and 
\item[(ii)]
$\hat{L}(\widehat{\sum a_{\sigma} {\sigma}} )=\mu(\sum a_\sigma \lambda(\sigma))
\hat{L}(\hat{g_0})$, \ for any $\widehat{\sum a_\sigma {\sigma}}\in 
\CC \Sigma/J$,
\end{enumerate}
if and only if the following conditions hold:
\begin{enumerate}
\item[(iii)] If $\sum_{j=1}^m a_{j}\alpha_j\in J$ and 
 $\omega, \beta\in *_{i=1}^n G_i^+$ are such that 
 $\{\omega \alpha_1 \beta,\ldots, \omega \alpha_m \beta\}\subset \Sigma$, then
\begin{equation}\label{eq1}
\sum_{j=1}^m a_j\hat{L}(\widehat{\omega\alpha_j \beta})=0.
\end{equation}
\item[(iv)]
$K_1\leq K_2$,  where the kernels
 $K_1$ and $K_2$ are associated, as in  \eqref{k1} and \eqref{k2},
  with the admissible  set
  $\Sigma $ and the map $L:\Sigma\to B(\cH) $ defined  by 
  $L(\sigma):= \hat{L}(\hat{\sigma})$, 
  \ $\sigma\in \Sigma$.
\end{enumerate}
\end{theorem}
  \begin{proof}
  Let $\hat{L}: \Sigma/J\to B(\cH)$ be such that the conditions (i) and (ii) hold. 
  For any $\alpha\in \Sigma$, we have
  $$L(\alpha)=\hat{L} (\hat \alpha )=\mu(\lambda(\alpha))\hat L(\hat g_0)= 
  \mu(\lambda(\alpha)) L(g_0).
  $$
  According to Theorem \ref{main-semigroup}, we deduce that $K_1\leq K_2$.
 Moreover, if $\sum_{j=1}^m a_{j}\alpha_j\in J$ and 
  $\omega, \beta\in \FF_n^+$ 
 are such that 
 $\{\omega \alpha_1 \beta,\ldots, \omega \alpha_m \beta\}\subset \Sigma$, then
 $$
 \sum_{j=1}^m a_j\hat{L}(\widehat{\omega\alpha_j \beta})=
 \mu(\lambda(\omega)) \mu\left(\sum_{j=1}^m a_j \lambda(\alpha_j)\right) 
 \mu(\lambda(\beta))
 \hat L(\hat g_0)=0.
 $$ 
 
 Conversely, assume the conditions (iii) and (iv) hold. As in the proof of Theorem 
 \ref{main-semigroup}, we define the representations $T_i$ of $G_i^+$ on $\cH$
  satisfying the condition \eqref{contrep}, the relation 
  $$
 T_i(g_i)L(\sigma)h=\begin{cases}L(g_i\sigma)h, & \text{ if } g_i\sigma\in \Sigma\\
 0, & \text{ otherwise},
 \end{cases}
 $$
 and 
  $T(\omega)|\cM^\perp=0$, $\omega\in *_{i=1}^n G_i^+\backslash \{g_0\}$,
   where $\cM$ is the closed span of all
  vectors $L(\sigma)h$ with $\sigma\in \Sigma$ and $h\in \cH$.
 Hence, we get $L(\omega)=T(\omega) L(g_0)$ for any $\omega\in \Sigma$.
 On the other hand, since $(\Sigma, J)$ is an admissible pair,  if 
  $\sum_{j=1}^m a_j\alpha_j\in J$, $\sigma\in \Sigma$, and $h\in \cH$, then 
$$
 \left(\sum_{j=1}^m a_jT({\alpha_j})\right) L(\sigma)h=
 \begin{cases}  \sum\limits_{j=1}^m a_j L(\alpha_j\sigma)h, & \text{ if }
  \{\alpha_1\sigma ,\ldots,  \alpha_m \sigma\}\subset \Sigma\\
 0, & \text{ otherwise}.
 \end{cases}
 $$
 Now, using relation \eqref{eq1}, we deduce 
 $\sum_{j=1}^m a_jT({\alpha_j})=0$  for 
 any $\sum_{j=1}^m a_j \alpha_j\in J$.
  Consider  
 $\mu: \cO(*_{i=1}^n G_i^+)\to B(\cH)$  to be the Poisson transform 
 associated with the representation
 $T$.
 Now, it is easy to see that, if $\omega\in \Sigma$, then
 $$
 \hat{L}(\hat{\omega})=L(\omega) =T_\omega L(g_0)= \mu(\lambda(\omega))
 \hat{L}(\hat{g_0}).
 $$
 The proof is complete.
\end{proof}

Let us remark that the Poisson transform $\mu$ of Theorem \ref{exten}
gives rise to a homomorphism 
 $\Psi:\CC *_{i=1}^n G_i^+/J\to B(\cH)$   by setting  
 $$
 \Psi(\widehat{\sum a_\sigma {\sigma}}):= \mu(\sum a_\sigma \lambda(\sigma)),
 $$
which  satisfies inequality \eqref{psihat}. The conditions (i) and (ii) of the 
 same theorem can be written as 
 $
 \hat L( f)= \Psi( f)\hat{L}(\hat{g_0})
 $
for any $ f\in 
\CC \Sigma/J$.
 Moreover, the homomorphism $\Psi$ give rise to a Hilbert module over 
$\CC *_{i=1}^n G_i^+/J$ in the natural way
$
f\cdot h:= \Psi(f)h$, \ $ f\in \CC *_{i=1}^n G_i^+/J$ ,  and  $h\in \cH$.

 Finally, we shoud mention that, as in Section \ref{Moment}, a vector-valued version 
of Theorem \ref{exten} can  easily be obtained.  This provides a characterization 
for the orbits  of Hilbert modules over $\CC *_{i=1}^n G_i^+/J$.

%



\end{document}